\documentclass{amsart}
\usepackage{amssymb}
\usepackage{verbatim}
\usepackage[all]{xy}

\newcommand{\short}{\mathcal{E}}
\newcommand{\cX}{\mathcal{X}}
\newcommand{\cY}{\mathcal{Y}}
\newcommand{\isoc}{\operatorname{Isoc^\dagger}}
\newcommand{\piun}{\pi_1^{\text{rig,un}}}
\newcommand{\aug}{\mathfrak{a}}
\newcommand{\Mod}{{\mathcal{M}od}}
\newcommand{\vvec}{{\mathcal{V}ec}}
\newcommand{\Rep}{{\mathbb{R}\textup{ep}}}
\newcommand{\crdot}{\odot}
\newcommand{\ctens}{\hat{\otimes}}
\newcommand{\utens}{\boxtimes}
\newcommand{\lie}{\operatorname{Lie}}

\newcommand{\hcheck}{\check{H}}

\newcommand{\del}{\nabla}
\renewcommand{\O}{\mathcal{O}}
\newcommand{\jdag}{j^\dagger}
\newcommand{\jdagY}{\jdag \O_{]Y[}}
\def\jdagOmY^#1{\jdag \Omega_{]Y[}^{#1}}
\def\dirlim_#1{\mathchoice{\underset{#1}{\varinjlim}}{\varinjlim_{#1}}{}{}}
\def\invlim_#1{\mathchoice{\underset{#1}{\varprojlim}}{\varprojlim_{#1}}{}{}}
\newcommand{\vv}{\mathcal{V}}
\newcommand{\shF}{\mathcal{F}}
\newcommand{\shG}{\mathcal{G}}
\newcommand{\Hom}{\operatorname{Hom}}
\newcommand{\spec}{\operatorname{Spec}}
\newcommand{\spf}{\operatorname{Spf}}
\newcommand{\Ext}{\operatorname{Ext}}
\newcommand{\pr}{^\prime}
\newcommand{\prpr}{^{\prime\prime}}
\newcommand{\rig}{{\textup{rig}}}
\newcommand{\hr}{H_{\rig}}
\newcommand{\dR}{{\textup{dR}}}
\newcommand{\one}{{\boldsymbol{1}}}
\newcommand{\isom}{\cong}
\newcommand{\hdr}{H_{\dR}}

\newcommand{\SP}{\operatorname{sp}}
\newcommand{\Ker}{\operatorname{Ker}}
\newcommand{\inject}{\hookrightarrow}
\newcommand{\loc}{{\textup{loc}}}
\newcommand{\col}{{\textup{Col}}}
\newcommand{\acol}{A_{\col}}
\newcommand{\alog}{A_{\textup{log}}}
\newcommand{\Olog}{\Omega_{\textup{log}}}
\newcommand{\aloc}{A_{\loc}}
\newcommand{\Ocol}{\Omega_{\col}}
\newcommand{\Oloc}{\Omega_{\loc}}

\newcommand{\aabs}{A_{\text{abs}}}
\newcommand{\id}{\operatorname{id}}

\newcommand{\F}{\mathbb{F}}
\newcommand{\C}{\mathbb{C}}

\newcommand{\Acol}{\O_{\col}}
\newcommand{\un}{\mathcal{U}n}
\newcommand{\ddb}{\operatorname{\bar{\partial}}}
\newcommand{\cU}{\mathcal{U}}
\newcommand{\ceck}{\v Cech}
\newcommand{\ninv}{\Psi}
\def\acoln_#1{A_{\col,#1}}
\def\Ocoln_#1{\Omega_{\col,#1}}
\def\iny(#1,#2){{#2}^{-\infty}#1}
\def\tu(#1,#2){{]#1[_{#2}}}
\def\hF(#1){H_{#1}^\otimes}

\newtheorem{theorem}{Theorem}[section]
\newtheorem{proposition}[theorem]{Proposition}
\newtheorem{lemma}[theorem]{Lemma}
\newtheorem{corollary}[theorem]{Corollary}
\theoremstyle{definition}
\newtheorem{definition}[theorem]{Definition}
\newtheorem{remark}[theorem]{Remark}

\numberwithin{equation}{section}

\begin{document}
\title[Coleman integration using the Tannakian formalism]{Coleman
integration using the Tannakian formalism}
\author{Amnon Besser}
\address{Department of Mathematics\\
Ben-Gurion University of the Negev\\
P.O.B. 653\\
Be'er-Sheva 84105\\
Israel
}
\maketitle

\section{Introduction}
\label{sec:intro}

The theory of Coleman integration has been developed in a
number of directions by Coleman~\cite{Col82,Col85,Col-de88} and later
by the 
author~\cite{Bes97} (a totally different approach was worked out by
Colmez~\cite{Colm96} and by Zarhin~\cite{Zar96}). It provides a way of
extending the ring of analytic
functions on certain rigid analytic domains in such a way that one can
often find an essentially unique primitive to a differential form. The
guiding principle has always been that $p$-adically one can integrate
differential forms ``locally'' but that due to the complete
discontinuity of the $p$-adic topology naive analytic continuation is
impossible and one must use instead ``analytic continuation along
Frobenius''.

The way in which Coleman does this analytic continuation,
using certain polynomials in an endomorphism which is a lift of
Frobenius, is admittedly rather strange. It has been our feeling that
this tends to put off people from exploring this theory in spite
of its many applications. This work grew out of an attempt to make the
theory more attractive to a wider mathematical audience by making the
procedure of continuation along Frobenius better understood. We later
realized that our methods makes it possible to develop the theory from
the start in a way which at least to us seems more elegant. In
addition, we are able to extend the theory of Coleman iterated
integrals, which was only developed for one dimensional spaces, to
spaces in arbitrary dimensions.

The idea behind the construction is very simple. An iterated Coleman
integral is the element $y_1$ in a solution of a system of
differential equations 
\begin{align*}
  d y_1 &= y_2\cdot \omega_1\\
  d y_2 &= y_3\cdot \omega_2\\
        &...\\
  d y_{n-1} &=  \omega_{n-1}\\
\intertext{and we may modify the last equation and add one more
  equation as follows:}
  d y_{n-1} &=  y_n\cdot \omega_{n-1}\\
  d y_{n} &=  0\;.
\end{align*}
We do not enter here into the nature of the differential forms
$\omega_i$.
Thus, a Coleman iterated integral is a solution of a unipotent system
of differential equations. Locally, this system can be solved and we
would like to analytically continue this solutions. Here, locally
means on a residue class $U_x$ - the set of points in our rigid space having
a common reduction $x$. According to a Philosophy of Deligne, the
collection of all paths between $x$ and another point $y$ along which
one can do analytic continuation is a principal homogeneous space for
a certain Tannakian fundamental group. There is an action of Frobenius
on this fundamental group and on the space of paths. Using work
of Chiarellotto~\cite{Chi96} we are able to show that there is a unique
Frobenius invariant path on this space and thus canonical
``analytic continuation along Frobenius'' has been achieved.

Let us get a bit more technical. Let $K$ be a discrete valuation field
with ring of integers $\vv$ and residue field $\kappa$ of
characteristic $p$. The basic setup for the theory is a
rigid triple $(X,Y,P)$ where $X$ is an open subscheme of the proper
$\kappa$-subscheme $Y$ which is in turn a closed subscheme of a
$p$-adic formal scheme $P$ smooth in a neighborhood of $X$. To $P$
one associates its generic fiber, a rigid analytic $K$-space, and in
it sits the tube (in the sense of Berthelot) $\tu(X,P)$ of points
reducing to points of $X$.

For the development of the theory we
restrict to the case where $\kappa$ is the algebraic closure of the
field $\F_p$ with $p$ elements. The main output of our constructions
is a certain ring $\acol(T)$, where $T$ is a shorthand for the triple
$(X,Y,P)$. This ring contains the ring of ``overconvergent'' rigid
analytic functions on $\tu(X,P)$. These rings form sheaves with
respect to a certain topology induced from the Zariski topology on
$X$. We can also define modules of differential forms $\Ocol^n(T)$ of
$n$-forms with coefficients in $\acol(T)$. One can form a \emph{Coleman de
Rham complex} from these differentials. The main result of the theory,
which is almost trivial once the machinery has been set up properly,
is that this complex of differentials is exact at the one-forms. This
corresponds to Coleman's theory where one construct functions by
successive integration in such a way that every differential form is
integrable. Note that in Coleman's theory one is dealing with curves so all
forms are automatically closed.

We can also interpret our Coleman functions as actual functions in the
style of Coleman. We get certain locally analytic functions on the
tubes $\tu(X,P)$. In Coleman's theory one gets a bit more since the
functions extend to a ``wide open space'' containing the tube. This
extension becomes complicated in higher dimensions so we have only
carried it out in dimension $1$. Indeed, we are able to show that our
theory is identical to Coleman's theory in the cases he considers.

This work begun while at the University of Durham and continued while
at the SFB 478 ``Geometrische Strukturen in der Mathematik'' at the
university of M\"unster. The author would like to thank both institutions. He
would also like to thank Bruno Chiarellotto J\"org Wildeshaus and
Amnon Yekutieli for some useful
suggestions and to Elmar Grosse Kl\"onne for reading parts of the
manuscript. Finally, he would like to thank the referee for a very
careful reading of the manuscript and for pointing out a number of
essential corrections.

Notation: 
$K$ is a non-archimedian field of characteristic $0$ with ring of
integers $\vv$ and residue
field $\kappa$ of characteristic $p$. We denote the category of finite
dimensional vector spaces over $K$ by
$\vvec_K$. All schemes are separated and of finite type over their
respective bases. If $X$ is a $\kappa$-scheme where $\kappa$ is
algebraic over the prime field $\F_p$, then
a \emph{Frobenius endomorphism} of $X$ is any morphism obtained in the
following way: Let $X\pr/\F_{p^r}$ be a scheme such that $X\isom
X\pr \otimes_{F_{p^r}} \kappa$. The morphism $\operatorname{Frob}^r$, where
$\operatorname{Frob}$ is the absolute Frobenius of $X\pr$, is an
$\F_{p^r}$-linear endomorphism of $X\pr$ and
$\operatorname{Frob}^r\otimes \id$, considered as a
$\kappa$-endomorphism of $X$ via the isomorphism above is a Frobenius
endomorphism of $X$.

\section{Unipotent isocrystals}
\label{sec:crystals}

In this section we review basic facts about connections and crystals,
in the rigid and overconvergent settings. Basic sources for these are
\cite{Ber96}, \cite{Chi-leS99} and \cite{Chi-leS99a}.

Let $V$ be a smooth rigid analytic $K$-space. If $E$ is a coherent sheaf
of $\O_V$-modules a connection on $E$ is a map 
\begin{equation*}
  \nabla: E \to E\otimes_{\O_V} \Omega_V^1
\end{equation*}
which satisfies the Leibniz formula
\begin{equation*}
  \nabla(fe)=e\otimes df + f\cdot \nabla(e)
\end{equation*}
for any $e\in E$ and $f\in \O_V$.
We will say that the pair
$(E,\nabla)$ is a connection on $V$. The connection $\nabla$ gives maps,
also denoted by $\nabla$, $\nabla:E\otimes_{\O_V} \Omega_V^i \to
E\otimes_{\O_V} \Omega_V^{i+1}$.
A connection $\nabla$ has a
curvature $\nabla^2\in \Hom_{\O_V}(E,E\otimes_{\O_V} \Omega_V^2)$ given by
\begin{equation}\label{curvature}
  \nabla^2 (e) = \nabla(\nabla(e))\;.
\end{equation}
A connection is called integrable if $\nabla^2=0$. In this case it is
well known that $E$ must be locally free.
When the connection on $E$ is integrable one obtains a de Rham complex
\begin{equation}\label{dRham}
  E \xrightarrow{\nabla} E\otimes_{\O_V} \Omega_V^1
  \xrightarrow{\nabla}  E\otimes_{\O_V} \Omega_V^2 \to \cdots\;. 
\end{equation}

A map of connections $m:(E,\nabla) \to (E\pr,\nabla\pr)$ is an $\O_V$-linear
map $m:E \to E\pr$ such that $\nabla\circ m = (m\otimes \id) \circ \nabla$.
This makes the collection of all connections on $V$ into a category.
The trivial connection is the connection $\one:=(\O_V,d)$. We will
also call a connection trivial if it is a direct sum of $\one$'s. A
connection is called unipotent if it is a successive extension of trivial
connections.

The following easy result will be the key to many of the constructions
to follow:
\begin{lemma}
  Let $(E,\nabla)$ be a connection with $E$ torsion free and let
  $f:E\to B$ be a map of $\O_V$-modules with $B$ torsion free as
  well. Let $A\subset E$ be the kernel of $f$. Construct by induction
  a sequence of subsheaves:
  \begin{equation*}
    A_0=A,\; A_{n+1}=\Ker \left[A_n \xrightarrow{\nabla}
    E\otimes_{\O_V} \Omega_V^1 \to  (E\otimes_{\O_V} \Omega_V^1) /
    (A_n\otimes_{\O_V} \Omega_V^1)\right] \;.
  \end{equation*}
  Let $A_\infty= \cap_{n=0}^\infty A_n$. Then 
  \begin{enumerate}
  \item $\nabla$ maps $A_\infty$ to $A_\infty \otimes_{\O_V}
    \Omega_V^1$.
  \item  The pair $\iny(A,\nabla):= (A_\infty,\nabla)$ is a
  connection.
\item  Any map of connections to $(E,\nabla)$ whose image lies in $A$
  factors through $\iny(A,\nabla)$.
  \end{enumerate}
\end{lemma}
\begin{proof}
The first and third assertions are clear. For the second assertion we
need to check that $A_\infty$ is coherent. The Leibniz formula implies
that the maps defining the $A_n$'s are $\O_V$-linear. Therefore, we see by
induction that $A_n$, being a kernel of an $\O_V$-map of coherent
sheaves and using~\cite[Proposition 9.4.3.2.ii]{BGR84}, is
coherent. Let $B_n:=E/A_n$. We show by induction on $n$ that
$B_n$ is also torsion free. For $n=0$ this is clear since $B_0\subset B$.
Suppose that $B_n$ is torsion free. From the short exact sequence
\begin{equation*}
  0\to A_n \to E \to B_n \to 0
\end{equation*}
and the fact that $\Omega_V^1$ is locally free since $V$ is smooth we find 
\begin{equation*}
  (E\otimes_{\O_V} \Omega_V^1)/(A_n\otimes_{\O_V} \Omega_V^1)
  = B_n \otimes_{\O_V} \Omega_V^1\;.
\end{equation*}
It follows that $A_n/A_{n+1}\subset B_n\otimes_{\O_V} \Omega_V^1$ and from
the induction hypothesis we get that $A_n/A_{n+1}$ is torsion free. The short
exact sequence
\begin{equation*}
  0 \to A_n/A_{n+1} \to B_{n+1} \to B_n \to 0
\end{equation*}
and the induction hypothesis now show that $B_{n+1}$ is an extension of
torsion free modules hence also torsion free.

In the course of the proof we showed that $A_n/A_{n+1}$ is torsion free for
each $n$. Consider now an open affinoid domain $U\subset V$. It
follows by rank consideration that the sequence $A_n|_U$ must
stabilize. Indeed, since both $A_n$ and $A_n/A_{n+1}$ are coherent, it
follows that their restrictions are associated to the
$\O_V(U)$-modules $A_n(U)$ and $A_n/A_{n+1}(U)$. By Noether's
normalization lemma \cite[Corollary 6.1.2.2]{BGR84} the algebra $\O_V(U)$
is finite over some Tate algebra $T_l$ and therefore
$\mathbb{L}:=\O_V(U)\otimes_{T_l}\operatorname{Frac}(T_l)$ is a
product of fields. By tensoring with $\mathbb{L}$ and considering the
dimensions over the various fields in the product it is clear that
$(A_n(U)/A_{n+1}(U))\otimes_{\O_V(U)}\mathbb{L}=0$ for sufficiently
large $n$, which implies, since $A_n(U)/A_{n+1}(U)$ is torsion free,
that $A_n(U)/A_{n+1}(U)=0$. It follows that $A_\infty|_U=A_n|_U$ for a
sufficiently large $n$ and is thus associated with the module
$A_\infty(U)$. This proves that $A_\infty$ is coherent and completes
the proof of the lemma.
\end{proof}

From the construction of $\iny(A,\nabla)$ the following results are clear.
\begin{lemma}\label{restriction}
  In the situation of the previous lemma suppose that $U\subset V$ is an
  open subspace. Then $\iny({(A|_U)},\nabla)=(\iny(A,\nabla))|_U$.
\end{lemma}
\begin{lemma}\label{restriction1}
  If, in the situation of the previous lemma, we have $y\in A(U)$ and
  $\nabla y=0$, then $y\in \iny(A,\nabla)(U)$.
\end{lemma}
We may apply the lemmas in particular to the kernel of the
curvature. We obtain a connection
\begin{equation}
  \label{eq:edel2}
  E^{\textup{int}}:= \iny(\Ker \nabla^2,\nabla)\;.
\end{equation}
The following are immediate consequences of the two lemmas.
\begin{corollary}\label{edel2cor1}
  Every map of connections from an integrable connection to
  $(E,\nabla)$ factors through $E^{\textup{int}}$.
\end{corollary}
\begin{corollary}\label{edel2cor2}
  If $U\subset V$ is an open subspace and $(E,\nabla)$ is a connection on
  $V$, then $E^{\textup{int}}|_U = (E|_U)^{\textup{int}}$.
\end{corollary}
\begin{corollary}\label{edel2cor3}
  If, in the situation of the previous corollary, we have $y\in E(U)$
  and $\nabla y=0$, then $y\in E^{\textup{int}}(U)$.
\end{corollary}
\begin{lemma}\label{untriv}
  If $V$ is an affinoid space and $E$ is a unipotent connection, then
  the underlying $\O_V$-module is free.
\end{lemma}
\begin{proof}
This is because on an affinoid space all extensions of $\O_V$ by
itself are trivial. Indeed, $\Ext^1(\O_V,\O_V)=H^1(V,\O_V)=0$.
\end{proof}
We now recall Berthelot's notions of tubes and strict neighborhoods
and some related terminology from~\cite{Bes98a}.
\begin{definition}[{\cite[Section~\ref{sec:rigidcomplex}]{Bes98a}}]
  A rigid triple is a triple $(X,Y,P)$ consisting of $P$, a formal $p$-adic
  $\vv$-scheme, $Y$ a closed $\kappa$-subscheme of $P$ which is proper
  over $\spec(\kappa)$ and $X$ an open $\kappa$-subscheme of $Y$ such
  that $P$ is smooth in a neighborhood of $X$.
\end{definition}

We will denote a rigid triple by a single letter, usually $T$. We let
$j$ denote the embedding of $X$ in $Y$.
When $T=(X,Y,P)$ is a rigid triple,
Berthelot defines the following spaces and notions. To $P$ corresponds
a rigid analytic space $P_K$, called the generic fiber of $P$,
together with a specialization map $\SP: P_K \to P$. To any locally
open subset $Z$ of $P$ (in particular to $X$ or
$Y$) corresponds a tube $\tu(Z,P)$, which is a rigid analytic subspace
of $P_K$ whose underlying set of points is the set $\SP^{-1}(Z)$. Let
$Z=Y-X$. An admissible open $U\subset \tu(Y,P)$ is called a strict
neighborhood of $\tu(X,P)$ if $(U,\tu(Z,P))$ is an admissible cover of
$\tu(Y,P)$.
Let $V$ be such a strict neighborhood.
Berthelot
defines a functor $\jdag$ from the category of sheaves on $V$
to itself by
\begin{equation*}
  \jdag(F) = \dirlim_{U} {j_U}_\ast F \; ,
\end{equation*}
where the direct limit is over all $U$ which are strict neighborhoods
of $\tu(X,P)$ in $\tu(Y,P)$ contained in $V$ and $j_U$ is the
canonical embedding.
We always have $(\jdag E)|_{\tu(X,P)}=E|_{\tu(X,P)}$.

A sheaf of $\jdag \O_{\tu(Y,P)}$-modules on $\tu(Y,P)$ will be called
a $\jdagY$-module for short. The notion of coherence for such modules
is recalled in \cite{Ber96} just before (2.1.9).
\begin{proposition}[{\cite[Proposition (2.1.10)]{Ber96}}]\label{jdagof}
  Any coherent $\jdagY$-module is obtained from a
  coherent $\O_V$-module for some strict neighborhood $V$ by applying
  $\jdag$ and any morphism between two such $\jdagY$-modules is $\jdag$ of a
  morphism over a possibly smaller strict neighborhood.
\end{proposition}
\begin{definition}[{\cite[\ref{rigid-triple-morphism}]{Bes98a}}]
  Rigid triples are made into a category in the following way: If
  $T\pr=(X\pr,Y\pr,P\pr)$ is another rigid triple and $f:X\to X\pr$ is
  a $\kappa$-morphism, a map of rigid spaces $F:U \to \tu(Y\pr,P\pr)$,
  where $U$ is a strict neighborhood of $\tu(X,P)$ in $\tu(Y,P)$, is
  said to be compatible with $f$ if its restriction to $\tu(X,P)$
  lands in $\tu(X\pr,P)$ and is compatible with $f$ via the
  specialization map. A morphism between $T$ and $T\pr$ consists then
  of a morphism $f$ together with a germ of a morphism $F$ compatible
  with $f$.
\end{definition}
It is proved in~\cite[Lemma~\ref{germs}]{Bes98a} that the inverse
image of a strict neighborhood of $\tu(X\pr,P\pr)$ in $\tu(Y\pr,P\pr)$
under a compatible morphism $F$ as above is a strict neighborhood of
$\tu(X,P)$ in $\tu(Y,P)$. This easily shows that morphisms of rigid
triples can be composed. This also shows that there is a pullback map
$f^\ast$ from $\jdag \O_{\tu(Y\pr,{})}$-modules to $\jdagY$-modules
along a morphism $f:T\to T\pr$.

A connection on a $\jdagY$-module $E$ is a map
\begin{equation*}
  \nabla: E \to E\otimes_{\O_{\tu(Y,P)}} \Omega_{\tu(Y,P)}^1
\end{equation*}
satisfying the Leibniz formula. Such a connection is called integrable if
its curvature, defined as in \eqref{curvature}, is $0$. The following
proposition is part (i) of \cite[Proposition (2.2.3)]{Ber96}.
\begin{proposition}
If $E$ is a $\jdagY$-module with an integrable connection $\nabla$, then there
exists a strict neighborhood $V$ of $\tu(X,P)$ and a connection
$(E_0,\nabla_0)$ on $V$ such that $(E,\nabla)=\jdag (E_0,\nabla_0)$.
\end{proposition}

Let $E$ be a coherent $\jdagY$-module with a connection $\nabla$ and
let $E\to B$ be a map to another $\jdagY$-module with kernel $A$.
suppose that there exists a strict neighborhood $V$ of $\tu(X,P)$,
a connection $(E_0,\nabla_0)$ on $V$ and a map of locally free $\O_V$-modules
$E_0 \to B_0$ such that $(E,\nabla)=\jdag
(E_0,\nabla_0)$ and $(E \to B)=\jdag(E_0 \to B_0)$. We have $A=\jdag A_0$ with 
$A_0=\Ker (E_0\to B_0)$.
\begin{proposition}
  The connection $\iny(A,\nabla):=\jdag \iny(A_0,\nabla_0)$ is
  independent of the choice of $V$, $E_0$, $\nabla_0$ and $B_0$.
\end{proposition}
\begin{proof}
  Immediate from Lemma~\ref{restriction}.
\end{proof}
We will abuse the standard terminology and call a $\jdagY$-module $E$
together with an integrable connection an isocrystal on the triple
$T$. An isocrystal is called overconvergent if the Taylor
expansion map gives an isomorphism on a strict neighborhood of the
diagonal. We denote the category of overconvergent isocrystals on $T$
by $\isoc(T)$. The category $\isoc(X,Y,P)$ is independent up to
equivalence of $Y$
and $P$ and therefore we may talk of overconvergent isocrystals on
$X$. We denote the category of overconvergent isocrystals on $X$ by
$\isoc(X,K)$ or simply by $\isoc(X)$ if $K$ is understood. If $f:X\pr
\to X$ is a $\kappa$-morphism there is a pullback map
from $\isoc(X)$ to $\isoc(X\pr)$ which is
defined in~\cite{Ber96} after Definition (2.3.6). If we have rigid
triples $T=(X,Y,P)$ and 
$T\pr=(X\pr,Y\pr,P\pr)$ and $f$ extends to a map of triples one can show
that this
pullback map can be realized as the pullback map along this morphism
of triples.

An overconvergent isocrystal over $\spec(\kappa)$ is nothing else than
a $K$-vector space, as can easily be seen by choosing the rigid triple
$(\spec(\kappa),\spec(\kappa),\spf(\vv))$. Given a point
$x:\spec(\kappa) \to X$ we can pullback along $x$ to obtain a functor
\begin{equation}
  \label{eq:fibfunct}
  \omega_x=x^\ast: \isoc(X) \to \vvec_K\;.
\end{equation}
To describe this explicitly it is useful 
to also consider the rigid triple 
\begin{equation}
  \label{tx}
  T_x:=(x,x,P)\;.
\end{equation}
Here of course, the scheme $x$ really means $\spec(\kappa)$ embedded
in $P$ via $x$. In this case the 
corresponding tube is the so called \emph{residue class} in Coleman's
terminology
\begin{equation}
  \label{eq:disc}
  U_x:= \tu(x,P)
\end{equation}
and the functor $\omega_x$ is obtained as
\begin{equation*}
  \omega_x(F,\nabla)=\{v\in F(U_x),\; \nabla(v)=0\}
\end{equation*}
If $f:X\pr \to X$ is a $\kappa$-morphism, $x\pr\in X\pr(\kappa)$ and
$f(x\pr)=x$, then there is a tautological natural isomorphism
\begin{equation}
  \label{tauisom}
  \omega_x \isom \omega_{x\pr} \circ f^\ast\;.
\end{equation}

An isocrystal on $X$ or on $T$ is called unipotent if it is a successive
extension of trivial connections. We denote the category of all
unipotent isocrystals on $X$ by $\un(X,K)$ or $\un(X)$ if $K$ is
understood. It follows from~\cite[Proposition 1.2.2]{Chi-leS99} that a
unipotent isocrystal is always
overconvergent. It follows easily from Proposition~\ref{jdagof} that
if $E$ is a unipotent isocrystal on $(X,Y,P)$ then there exists a
strict neighborhood $U$ of $\tu(X,P)$ and an integrable unipotent
connection on $U$ whose $\jdag$ gives $E$.
\begin{proposition}[{\cite[Proposition I.2.3.2]{Chi-leS99}}]\label{unitannak}
  The category $\un(X)$ is closed under extensions, subobjects,
  quotients, tensor products internal homs duals and inverse
  images. A point $x\in X(\kappa)$ supplies a fiber functor $\omega_x$
  making $\un(X)$ a neutral Tannakian category.
\end{proposition}

The cohomology of an overconvergent isocrystal $E$ on $T=(X,Y,P)$ is
given by $H^i(\tu(Y,P),\jdag (E\to E\otimes \Omega^1 \to\cdots))$, where the
complex in brackets is the de Rham complex of $E$. This cohomology
depends only on $X$ hence can be denoted by $H^i(X,E)$.

From now until the end of this section we assume that $X$ is smooth.
\begin{lemma}
  Let $U$ be an open dense subset of $X$ and let $E$ be a unipotent
  isocrystal on $X$. Then  the map 
  $H^0(X,E) \to H^0(U,E)$ is an isomorphism.  
\end{lemma}
\begin{proof}
This follows if we show
that if $D$ is the complement of $U$ then $H_D^1(X,E)=0$. This is true
for the trivial isocrystal by purity~\cite[Corollaire 5.7]{Ber97} and
then follows easily for all unipotent isocrystals by induction.
\end{proof}
\begin{corollary}\label{Hom-ext}
  In the situation of the previous lemma suppose that $F$ is another
  unipotent isocrystal. Then we have
  $\Hom(E,F) \isom \Hom(E|_U,F|_U)$.
\end{corollary}
\begin{proof}
This follows since $\Hom(E,F)=H^0(E\otimes F^\ast)$ and $E\otimes F^\ast$
is  unipotent by~Proposition~\ref{unitannak}.
\end{proof}
\begin{proposition}\label{crist-ext}
Let $U$ be an open
subset of $X$. Let $E$ be a unipotent isocrystal on $X$ and let $G$ be a
subcrystal of $E|_U$. Then $G$ extends to $X$, i.e., there exists a
subcrystal $G\pr$ of $E$ such that $G\pr|_U =G$.
\end{proposition}
\begin{proof}
We may assume that $U$ is dense in $X$. Otherwise, since $X$ is smooth
we can extend $G$ to the connected
components with non-empty intersection with it and extend it as zero
in the other components. We use induction on the rank of $E$.
We have a short exact sequence 
\begin{equation}
  \label{starstar}
  0\to \one \to E \to E\pr \to 0
\end{equation}
with $E\pr$ unipotent as well. Let $G\prpr$ be the image in $E\pr|_U$
of $G$. By the induction hypothesis it is the restriction to $U$ of a
subcrystal $E\prpr$ of $E\pr$. Pulling back \eqref{starstar} via
$E\prpr \to E\pr$ we see that we may assume that $G\to E\pr|_U$ is
surjective. Now by rank considerations, either $G=E|_U$, in which case
the result is clear, or $G\to E\pr|_U$ is an isomorphism. In this last
case we find a splitting
$E\pr \to E$ over $U$ and by Corollary~\ref{Hom-ext} this extends to $X$
giving us our required subcrystal.
\end{proof}

\section{Unipotent groups}\label{sec:unip}
As explained in the previous section the category $\un(X)$ of unipotent
isocrystals on a smooth scheme $X$ is a Tannakian category and to every
point $x\in X(\kappa)$ corresponds a fiber functor
$\omega_x$. At this point we do assume we have such a point. We can
always achieve that by making an extension of scalars. From the theory of
Tannakian categories~\cite{Del-Mil82} it follows
that to $\un(X)$ corresponds an affine group scheme $G=\piun(X,x)$.
Let $F$ be a Frobenius endomorphism on
$X$ fixing $x$ (such an automorphism always exists). It
induces an automorphism of $G$, denoted $\phi$. If $y\in X(\kappa)$ is
another point of $X$, then there is a principal $G$-homogeneous
space $P_{x,y}$ representing the functor
$\operatorname{Isom}(\omega_x,\omega_y)$. In particular, an element of
$P_{x,y}(K)$ consists of a family of isomorphisms
$\lambda_E:\omega_x(E) \to \omega_y(E)$ indexed by the objects of
$\un(X)$ and satisfying the following properties:
\begin{enumerate}
\item The map $\lambda_\one$ is the identity on $K$.
\item For every two objects $M$ and $N$ of $\un(X)$ we have
  $\lambda_{M\otimes N} = \lambda_M \otimes\lambda_N$.
\item For every map $\alpha:M\to N$ we have $\lambda_N \circ
  \omega_x(\alpha)= \omega_y(\alpha) \circ \lambda_M$.
\end{enumerate}
There is an obvious composition map
\begin{equation}
  \label{composmap}
  P_{x,y} \times P_{y,z} \to P_{x,z}\;.
\end{equation}
If $f:X\to Y$ is a morphism with $f(x\pr)=x$ and $f(y\pr)=y$, there is
an induced map $f_\ast:P_{x\pr,y\pr} \to P_{x,y}$. On $K$-points this
map sends the collection of maps $\lambda_E$ to the collection of maps
\[
  f_\ast(\lambda)_E:=\lambda_{f^\ast E}:
  \omega_x(E)=\omega_{x\pr}(f^\ast E) \to
  \omega_{y\pr}(f^\ast E)=\omega_y(E)\;.
\]
Suppose our Frobenius endomorphism fixes both $x$ and $y$.
The action of Frobenius then induces an
automorphism of $P_{x,y}$, denoted $\varphi$, compatible with
the $G$-action in the sense that $\varphi(ga)=\phi(g)\varphi(a)$
for $g\in G$ and $x\in P$. The goal of this section is to prove the
following theorem.
\begin{theorem}\label{Tanakthm}
  The map $G\to G$ given by
  $g\mapsto g^{-1} \phi(g)$ is an isomorphism.
\end{theorem}
\begin{corollary}\label{prevcor}
  For any $x,y\in X(\kappa)$ there is a unique $a_{x,y}\in P_{x,y}(K)$ fixed by
  the automorphism $\varphi$ induced by any Frobenius endomorphism $F$
  fixing both $x$ and $y$. Furthermore, The composition map
  \eqref{composmap} sends $(a_{x,y},a_{y,z})$ to $a_{x,z}$.
\end{corollary}
\begin{proof}
Any two Frobenius endomorphisms have a common power so that it is
sufficient to prove the result for a single Frobenius endomorphism.
We first show uniqueness. In fact we show this for a point of $P_{x,y}$
in an arbitrary extension field $L$. Suppose $a$ and $ga$ are both
fixed by $\varphi$ with $g\in G(L)$.
Then we have 
\[
  ga=\varphi(ga) = \phi(g)\varphi(a)=\phi(g) a\;.
\]
Therefore, $g=\phi(g)$ so $g^{-1} \phi(g) = e = e^{-1} \phi(e)$ so the
theorem implies $g=e$.
We now show existence. For some possibly huge field $L/K$, which we
may assume Galois by further extension, we have a
point $a_0\in P_{x,y}(L)$. For the
point $a=ga_0$ to be a fixed point we need
\[
  ga_0 = \varphi(ga_0)=\phi(g)\varphi(a_0)\;,
\]
or $g^{-1}\phi(g) \varphi(a_0)=a_0$. Let $h\in G(L)$ be the unique element
such that $h\varphi(a_0)=a_0$. By the theorem we may solve 
$g^{-1}\phi(g) =h$ and obtain our fixed point $a=a_{x,y}$. Since the
fixed point $a$ is unique, we have $\sigma(a)=a$ for
every Galois automorphism, which shows that $a$ is already defined over
$K$. The behavior of these points under composition is again clear
from uniqueness.
\end{proof}
We can reformulate the above result in much more elementary terms
using the residue classes of \eqref{eq:disc}.
\begin{corollary}\label{unip-res}
  For every unipotent overconvergent isocrystal $E$ on a rigid
  triple$(X,Z,P)$, for every two
  points $x,y\in X(\kappa)$ and for every $v_x\in E(U_x)$ with
  $\nabla(v_x)=0$ there is a unique way to associate $v_y \in E(U_y)$
  with $\nabla(v_y)=0$ such that the following properties are
  satisfied:
  \begin{enumerate}
  \item\label{unr1} The association $v_x \mapsto v_y$ is $K$-linear.
  \item\label{unr2} For the trivial isocrystal it sends $1$ to $1$.
  \item\label{unr3} It is functorial in the following sense: If
    $f:E\to E\pr$ is a 
    morphism of isocrystals then $f(v_y)$ is associated with $f(v_x)$.
  \item\label{unr4} It is compatible with tensor products in the
    following sense:
    If $v_x\pr$ is a horizontal section of the isocrystal $E\pr$ on $U_x$
    associated with the horizontal section $v_y\pr\in E\pr(U_y)$, then
    the section $v_x\otimes v_x\pr \in E\otimes E\pr(U_x)$ is
    associated with $v_y\otimes v_y\pr \in E\otimes E\pr(U_y)$.
  \item\label{unr5} It is compatible with Frobenius in the following
    sense: Let $F:X\to X$ be a Frobenius endomorphism that fixes the
    points $x$ and $y$. Then we have the commutative diagram
    \begin{equation*}
      \xymatrix{
        \omega_x(E)\ar[r]\ar[d] &\omega_y(E)\ar[d] \\
        \omega_{x}(F^\ast(E))\ar[r] &\omega_{y}(F^\ast(E))
        }
    \end{equation*}
    with the vertical isomorphisms \eqref{tauisom}.
  \end{enumerate}
  In addition this association satisfies that if $v_y$ is associated
  with $v_x$ and $v_z$ is associated with $v_y$, then $v_z$ is
  associated with $v_x$.
\end{corollary}
\begin{proof}
This is just a reformulation of corollary~\ref{prevcor}.
\end{proof}

The last property above is valid for any morphism, not just a
Frobenius endomorphism, and in still greater generality.
\begin{proposition}\label{pullback}
  Let $f:X\to Y$ be a morphism and suppose
  $f(x\pr)=x$ and $f(y\pr)=y$. Let $E\in \un(Y)$. Then we have the
  commutative diagram
  \begin{equation*}
    \xymatrix{
      \omega_x(E)\ar[r]\ar[d] &\omega_y(E)\ar[d] \\
      \omega_{x\pr}(f^\ast(E))\ar[r] &\omega_{y\pr}(f^\ast(E))
      }
  \end{equation*}
\end{proposition}
\begin{proof}
In the diagram above, the horizontal maps are induced by the elements
$a_{x,y}\in P_{x,y}$ and $a_{x\pr,y\pr}\in P_{x\pr,y\pr}$ in the top
and bottom row respectively. The vertical maps are given in \eqref{tauisom}
For appropriate choice of Frobenius endomorphisms $F_X$ and $F_Y$ on
$X$ and $Y$ respectively, inducing maps $\varphi_X$ and $\varphi_Y$ on
the corresponding principal spaces, we have $F_Y\circ f=f\circ
F_X$. It follows from this that $f_\ast(a_{x\pr,y\pr})$ is $\varphi_Y$
invariant and therefore must be equal to $a_{x,y}$. Spelling this out
gives the result.
\end{proof}
\begin{definition}\label{continuation}
  We call the association $v_x \mapsto v_y$ analytic continuation
  along Frobenius from $x$ to $y$.
\end{definition}

To begin the proof of the theorem we need to recall the work of
Chiarellotto~\cite{Chi96}. Chiarellotto considers the completed
enveloping algebra $U=\hat{\mathcal{U}}(\lie G)$
of the Lie algebra of the group $G$ 
together with its augmentation ideal $\aug$. We summarize the results
of loc.\ cit.\ that we need in the following theorem.
\begin{theorem}\label{Chiarthm}
  \begin{enumerate}
  \item\label{chtm1}
    The $K$-algebra $U$ is complete with respect to the $\aug$-adic
    topology.
  \item\label{chtm2}
    The quotients $U_n=U/\aug^n$ are finite dimensional over $K$.
  \item\label{chtm3}
    The category $\Rep_K(G)$ of algebraic representations of $G$
    over $K$ is equivalent 
    to the category $\Mod_U$ of $U$-modules which are finite dimensional
    $K$-vector spaces and for which the action of $U$ is continuous
    when they are provided with the discrete topology. In particular,
    the category $\Mod_U$ is Tannakian and the functor $\omega: \Mod_U
    \to \vvec_K$ sending a module to its underlying vector space is a
    fiber functor.
  \item\label{chtm4} There is an automorphism $\phi:U\to U$ of the augmented
    algebra $U$ such that under the
    equivalence $\Mod_U \isom \Rep_K(G)\isom \un(X)$ the operations of
    twisting a $U$-module by $\phi$, twisting a $G$-representation by
    $\phi$ and the pullback $F^\ast$ via $F$ correspond.
  \item\label{chtm5}
    the $\phi$-module $\aug^n/\aug^{n+1}$ is mixed with negative
    weights, i.e., $\phi$ has on it eigenvalues which are Weil
    numbers with strictly  negative weights.
  \end{enumerate}
\end{theorem}
\begin{proof}
Part~\ref{chtm1} is \cite[Lemma II.2.4]{Chi96}. \ref{chtm2}.\  and
\ref{chtm5}.\  
both follow from the fact, proved in the course of proving
Proposition~II.3.3.1 there, that there is a surjective morphism, Frobenius
equivariant from its construction, $\hr^1(X/K)^{\otimes (-n)} \to
\aug^n/\aug^{n+1}$, and it is well known that $\hr^1(X/K)$ is mixed
with positive weights.

Part \ref{chtm3} is stated in \cite[II, \S 2]{Chi96} and proved
in the appendix to~\cite{Chi-leS99a}.
Finally, \ref{chtm4}.\ is discussed in section II.3 of \cite{Chi96}.
\end{proof}
We denote by $\pi_{n,m}: U_n \to U_m$ the obvious projection.
Notice that the completeness with respect to the $\aug$-adic topology
means that $U=\invlim_n U_n$.
We denote by $U\ctens U$ the completed tensor product which can be
defined by 
\begin{equation}\label{ctens}
  U\ctens U:=\invlim_n U_n \otimes U_n\;.
\end{equation}
The next result is fairly obvious and probably well known, giving $U$
part of the structure of a coalgebra.
\begin{lemma}\label{coprod}
  There exists a canonical algebra homomorphism $\Delta: U \to U\ctens
  U$ such that the
  tensor product in $\Mod_U$, which we denote by $\utens$, is given in
  the following way: If $M$ and
  $N$ are two $U$-modules, then $M\utens N$ has the tensor
  product over $K$ as an underlying vector space and the $U$-action is
  the composition of $\Delta$ with the obvious $U\ctens U$-action.
\end{lemma}
\begin{proof}
That $M\utens N$ has the underlying $K$-vector space structure
$M\otimes_K N$ is clear from the fact that $\omega$ is a tensor
functor.
By Theorem~\ref{Chiarthm}.\ref{chtm2} we have $U_n \in \Mod_U$, hence
$U_n \utens U_n\in \Mod_U$. Write the action of $U$ on objects of
$\Mod_U$ as $(u,m)\mapsto u\crdot m$. Let
$u\in U$. Then we define a sequence $\Delta_n(u):= u\crdot(1\otimes
1)\in U_n\otimes U_n$,
where $1\otimes 1 \in U_n\otimes U_n$. It is easy to see that the
$\Delta_n(u)$ define an element $\Delta(u)\in U\ctens U$ and that
$\Delta$ so defined is $K$-linear. Suppose now that $M$ and $N$ belong
to $\Mod_U$. The $U$-action factors through some $U_n$. For each $x\in
M$ and $y\in N$ there are maps of $U$-modules, $u \mapsto ux$ and
$u\mapsto uy$, from $U_n$ to $M$ and $N$ respectively. Their tensor
product sends $1\otimes 1$ to $x\otimes y$ hence $\Delta_n(u)=u\crdot
(1\otimes 1)$ to $u\crdot (x\otimes y)$ so $u\crdot (x\otimes y)=
\Delta_n(u) (x\otimes y)=\Delta(u) (x\otimes y)$. This is the last
assertion of the lemma from which the fact that $\Delta$ is a ring
homomorphism follows easily. Indeed, we have
\begin{equation*}
  \Delta_n(uv)=(uv)\crdot (1\otimes 1) = u\crdot (v\crdot (1\otimes
  1)) = \Delta_n(u) \Delta_n(v) (1\otimes 1)=\Delta_n(u) \Delta_n(v)\;.
\end{equation*}
\end{proof}

Now let $L$ be a commutative $K$-algebra. We have $U_L:=U\ctens L=\invlim_n
U_n \otimes_K L$. The maps $\epsilon$ and $\Delta$ extend to
$\epsilon_L:U_L \to L$ and $\Delta_L:U_L \to U_L\ctens_L U_L$.
\begin{proposition}
  Let $G\pr$ be the functor from commutative $K$-algebras to groups
  given by
  \begin{equation*}
    G\pr (L) =\{u\in U_L|\; \epsilon_L(u)=1,\; \Delta_L(u)=u\otimes u\}\;,
  \end{equation*}
  which is a group with the multiplication induced by the algebra
  multiplication on $U_L$. Then there is a natural isomorphism
  $G\isom G\pr$
\end{proposition}
\begin{proof}
By Theorem~\ref{Chiarthm}.\ref{chtm3} the group $G$ is the group
corresponding to the Tannakian category $\Mod_U$ together with the
fiber functor $\omega$ sending a module to its underlying vector
space.
We recall that an element of $G(L)$ consists of a family of maps
$(\lambda_M)$, where $M$ runs over the objects of $\Mod_U$. For each $M$,
$\lambda_M: \omega(M)\otimes_K L \to \omega(M)\otimes_K L$ is an
automorphism. Such a
family should satisfy the following properties:
\begin{enumerate}
\item The map $\lambda_\one$ is the identity on $L$, where $\one$ is
  the unit object of $\Mod_U$, i.e, the $U$-module $U/\aug$.
\item\label{prop2} For any two objects $M$ and $N$ we have
  $\lambda_{M\utens N} =
  \lambda_M \otimes_L\lambda_N$.
\item For every map $\alpha:M\to N$ we have $\lambda_N \circ
  (\omega(\alpha)\otimes 1)= (\omega(\alpha)\otimes 1) \circ \lambda_M$.
\end{enumerate}
It is clear that $G\pr (L)$ is a group.
Let $M\in \Mod_U$. Clearly $U_L$ acts on $M\otimes_K L$.
An element $u\in G\pr (L)$ therefore defines an automorphism
$\lambda_M$ of
$M\otimes_K L$. It is immediately checked that the family $(\lambda_M)$
defines an element of $G(L)$. This gives a map $G\pr (L) \to
G(L)$. The inverse map is given as follows: Let $(\lambda_M)$ be an
element of $G(L)$. We define $u_n=\lambda_{U_n}(1) \in U_n\otimes_K
L$. The elements $u_n$ are compatible under the maps $\pi_{n,m}\otimes
1$ since the identity elements are. They thus define an element $u\in
U_L$. The inverse map will map $(\lambda_M)$ to $u$ and this map is
well defined once we show that $u\in G\pr (L)$. That $\epsilon_L(u)=1$
follows immediately from $\lambda_\one = 1$.
We next show that for any $M\in \Mod_U$ the automorphism $\lambda_M$
is given by multiplication by $u$. Suppose $M$ is in fact a
$U_n$-module. Let $x\in M$ and define a map $\alpha_x:U_n \to M$ by
$\alpha_x(v)=vx$. This is clearly a map of $U$-modules. Then
\begin{equation*}
  \lambda_M(x\otimes 1)=\lambda_M \circ (\alpha_x\otimes 1) (1) =
  (\alpha_x \otimes 1) \circ \lambda_{U_n}(1) =(\alpha_x \otimes 1) u_n =
  u_n \cdot (x\otimes 1) = u \cdot (x\otimes 1)\;.
\end{equation*}
By $L$-linearity this extends to give the claim. This shows
that the composition $G(L) \to G\pr (L) \to G(L)$ is the identity. It
remains to check that $\Delta(u) = u\otimes u$ and for that it
suffices to check that $\Delta_n(u_n)=u_n\otimes u_n$. But this is now
clear from property~\ref{prop2} of the family $\lambda_M$ and from the
construction of $\Delta_n$.
\end{proof}
We henceforth identify $G$ with $G\pr$. By
Theorem~\ref{Chiarthm}.\ref{chtm4} twisting the $U$-action by $\phi$ is a
tensor functor, from which it follows via Lemma~\ref{coprod} that
$\Delta\circ \phi=\phi \ctens \phi \circ \Delta$.
This immediately shows that the obvious action of $\phi$ on $U$
induces an action of $\phi$ on $G\pr$. In fact, via the identification
of $G$ with $G\pr$ this action is exactly the action of $\phi$ on $G$,
as can be straightforwardly checked after noting that this last action
is given by sending the collection $(\lambda_M)$ to the collection
$(\mu_M)$ where $\mu_M=\lambda_{F^\ast M}$ ($F^\ast M$  has the same
underlying vector space as $M$) and using
Theorem~\ref{Chiarthm}.\ref{chtm4} again.
\begin{proof}[Proof of Theorem~\ref{Tanakthm}]
We need to show that $g\mapsto g^{-1} \phi(g)$ induces a bijection on
$L$-points for any commutative $K$-algebra $L$.
We first claim that for any $a\in U_L$ with $\epsilon_L(a)=1$ there is
a unique $x\in U_L$ with $\epsilon_L(x)=1$ such that $\phi(x)=xa$. We
show by induction the same for $x_n\in U_n\otimes L$. Then by
uniqueness these solutions glue to give the unique $x$. For $n=1$
there is nothing to prove. Suppose we already found $x_{n-1}$. Let
$B=\{x\in U_n\otimes L|\; \pi_{n,n-1}(x)=x_{n-1}\}$. This is an affine
space for $C\otimes_K L$ where $C:=\aug^{n-1}/\aug^n$. We have a map
$T:B\to C\otimes_K L$ given by
$T(x)=\phi(x)-xa$. Since $\epsilon_L(a)=1$ we find for $y\in
C\otimes_K L$ that
$T(y+x)=(S\otimes 1)(y)+T(x)$ with $S:C\to C$ given by $S(y)=\phi(y)-y$. It
follows from Theorem~\ref{Chiarthm}.\ref{chtm5} that $S$ is
invertible. This immediately implies that $T$ is invertible and 
in particular there is a unique solution $x_n$ to $T(x_n)=0$.

It remains to prove that if $a\in G(L)$, then the unique $x$ so
constructed is also in $G(L)$. To do this we repeat the above argument
with $U_L$ replaced by $U_L\ctens U_L$, $\phi$ replaced by
$\phi\ctens \phi$ and $\epsilon_L$ replaced by $\epsilon_L \cdot
\epsilon_L$. The same argument works by using \eqref{ctens} and the
fact that the successive quotients $U_n\otimes U_n/U_{n-1}\otimes
U_{n-1}$ still have only negative weights. We need to
prove that $\Delta(x)=x\otimes x$. We have
\begin{equation*}
  \Delta(x) (a\otimes a) = \Delta(x)
  \Delta(a)=\Delta(xa)=\Delta(\phi(x))= (\phi\ctens \phi) \Delta(x)\;.
\end{equation*}
On the other hand
\begin{equation*}
  (x\otimes x) (a\otimes a) = (xa)\otimes (xa)=\phi(x)\otimes \phi(x)
  =(\phi\ctens \phi)(x\otimes x)\;.
\end{equation*}
We see that $\Delta(x)$ and $x\otimes x$ are both solutions to the
same equation which,
since $(\epsilon_L \cdot \epsilon_L)(a\otimes a)=1$ has a unique
solution, and they are therefore equal.
\end{proof}

\section{Coleman functions}
\label{sec:functions}
Starting from this section we assume that $\kappa$ is the algebraic
closure of $\F_p$. This restriction is imposed by wanting $\kappa$ to
be algebraic over $\F_p$ and also algebraically closed. The first
requirement seems essential because we are using Frobenius
endomorphisms. The second restriction is not essential but it
makes some things more pleasant.

Let $T=(X,Y,P)$ be a rigid triple as in section~\ref{sec:crystals}. As we
did in section~\ref{sec:unip} we denote, for a point $x\in X$ (a point
will always mean a closed point), the 
corresponding tube $\tu(x,P)$ by $U_x$. For a locally free
$\jdagY$-module $\shF$ we define $\shF(T)=\Gamma(\tu(Y,P),\shF)$. We
also denote $A(T)=\Gamma(\tu(Y,P),\jdagY)$ and
$\Omega^i(T)=\Gamma(\tu(Y,P),\jdagOmY^i)$.

\begin{definition}
  Let $\shF$ be a locally free $\jdagY$-module. 
  The category $\aabs(T,\shF)$ of \emph{abstract Coleman functions} on $T$
  with values in 
  $\shF$ is defined as follows: Its objects are triples $(M,s,y)$ where
  \begin{itemize}
  \item $M=(M,\del)$ is a unipotent isocrystal on $T$.
  \item $s\in \Hom(M,\shF)$.
  \item $y$ is a collection of sections, $\{y_x\in M(U_x),\; x\in
    X\}$, with
    $\del(y_x)=0$, which correspond to each other via ``analytic
    continuation along Frobenius'' as in Definition~\ref{continuation}.
  \end{itemize}
  A homomorphism $f$ between $(M,s,y)$ and $(M\pr,s\pr,y\pr)$ is a
  morphism of isocrystals $f:M\to M\pr$ such that $f^\ast(s\pr)=s$ and
  $f(y_x)=y_x\pr$ for any $x\in X$.
\end{definition}
Note that the condition $f(y_x)=y_x\pr$ need only be
checked at one point and then holds for all other points using
Corollary~\ref{unip-res}.\ref{unr3}. 
The direct sum of two abstract Coleman functions is given by
\begin{equation*}
  (M_1,s_1,y_1)+
  (M_2,s_2,y_2)= 
  (M_1\bigoplus M_2,s_1+s_2,y_1\oplus y_2)\;.
\end{equation*}
\begin{definition}
Given two locally free $\jdagY$-modules  $\shF$ and $\shG$ we define a
tensor product 
functor $\otimes: \aabs(U,\shF) \times \aabs(U,\shG) \to
\aabs(U,\shF\otimes \shG)$ by the formula
\begin{equation*}
  (M_1,s_1,y_1)\otimes
  (M_2,s_2,y_2):= 
  (M_1\bigotimes M_2,s_1\otimes s_2,y_1\otimes y_2)\;.
\end{equation*}
\end{definition}
\begin{definition}
  Let $T$ and $\shF$ be as above. The collection $\acol(T,\shF)$ of
  \emph{Coleman functions} 
  on $T$ with values in $\shF$ is the collection of connected
  components of the category $\aabs(T,\shF)$. The Coleman function
  corresponding to the triple $(E,s,y)\in\aabs(T,\shF)$ we denote by
  $[E,s,y]$. In particular we define Coleman functions and forms and
  Coleman forms with values in a $\jdagY$-module,
  \begin{align*}
    \acol(T):&=\acol(T,\jdagY)\;,\\
    \Ocol^i(T):&= \acol(T,\jdagOmY^i)\;,\\
    \Ocol^i(T,\shF):&= \acol(T,\shF\otimes_{\jdagY} \jdagOmY^i)\;.
  \end{align*}
\end{definition}
It will follow from Proposition~\ref{pullinj} below that the
collection of Coleman functions is in fact a set.
The following proposition is formal and is left for the reader.
\begin{proposition}
  The direct sum of abstract Coleman functions gives $\acol(T,\shF)$
  the structure of an abelian group and it is in fact a $K$-vector
  space with multiplication by a scalar $\alpha$ given by multiplying
  (say) the third component with $\alpha$. If $\shG$ is another
  locally free $\jdagY$-module then the tensor product induces a
  bilinear map $\acol(T,\shF)\times \acol(T,\shG) \to
  \acol(T,\shF\otimes \shG)$. In particular, the set $\acol(T)$ is a
  commutative $K$-algebra.
\end{proposition}
\begin{definition}\label{changesheaf}
  Let $f:\shF \to \shG$ be a map of $\jdagY$-modules. Then we obtain a
  functor $f_\ast: \aabs(T,\shF) \to  \aabs(T,\shG)$ via
  $f_\ast(M,s,y)= (M,f(s),y)$, and therefore a map $f_\ast: \acol(T,\shF)
  \to \acol(T,\shG)$.
\end{definition}

Recall that if $(F,\nabla_F)$ and $(M,\nabla_M)$ are two connections,
then the induced connection on  $\Hom(M,F)$ is given by
\begin{equation*}
  (\nabla_{\Hom(M,F)}f)(m)= \nabla_F(f(m)) - f\otimes \id_{\Omega^1}
  \nabla_M(m)\;.
\end{equation*}
\begin{definition}
Let $(F,\nabla_F)$ be a connection and let $(M,\nabla_M)$ be a
unipotent isocrystal on $T$. Using the de Rham differentials of $\Hom(M,F)$
we obtain the de Rham differentials on Coleman forms with values in
$F$, $\nabla_F:\Ocol^i(T,F) \to \Ocol^{i+1}(T,F)$ given by
$\nabla_F[M,s,y]=[M,\nabla_{\Hom(M,F)} s,y]$. If $\nabla_F$ is
integrable, then the de Rham differentials on Coleman forms give a
complex,
\begin{equation*}
  \acol(T,F) \to \Ocol^1(T,F)\to \Ocol^2(T,F)\to \cdots,
\end{equation*}
called the \emph{Coleman de Rham complex} of $(F,\nabla_F)$. In particular,
for the trivial connection we obtain the Coleman de Rham complex of $T$.
\end{definition}
\begin{definition}
Let $f:T\pr \to T$ be a map of triples.
Then there is a pullback functor $f^\ast:\aabs(T,\shF) \to
\aabs(T\pr,f^\ast \shF)$ given by $f^\ast(E,s,y)=(f^\ast E,f^\ast
s,f^\ast y)$. This is well defined since Proposition~\ref{pullback}
shows that the pullbacks of sections corresponding under 
analytic continuation on $T$ give
sections on $T\pr$ with the same property. This functor induces a $K$-linear map
$f^\ast:\acol(T,\shF) \to \acol(T\pr,f^\ast \shF)$. Using
Definition~\ref{changesheaf} we
immediately obtain a ring homomorphism $f^\ast:\acol(T) \to \acol(T\pr)$
and maps $f^\ast \Ocol^i(T) \to \Ocol^i(T)$ compatible with the
differentials of the Coleman de Rham complex.
\end{definition}

If we wish to interpret Coleman functions as actual functions, in the
style of Coleman, we can do the following.

\begin{lemma}\label{indep-of-rep}
  Suppose $X=x$ is a point. Let $[M,s,y]\in \acol(T,\shF)$. Then the
  section 
  \begin{equation*}
    \theta([M,s,y]):= s(y_x)\in \shF(U_x)
  \end{equation*}
  depends only on the Coleman function and not on the particular
  representing abstract Coleman function.
\end{lemma}
\begin{definition}
  The space of locally analytic functions on $T$ with values in $\shF$
  is defined to be the product
  \begin{equation*}
    \aloc(T,\shF):=\prod_{x\in X} \shF(U_x)\;.
  \end{equation*}
\end{definition}
Note that $\aloc(T,\jdagY)$ is naturally a ring.
\begin{definition}
  The map $\theta: \acol(T,\shF) \to \aloc(T,\shF)$ is defined by 
  \begin{equation*}
    f\to \prod_{x\in X} \theta(x^\ast f)\;.
  \end{equation*}
\end{definition}
\begin{proposition}
  The map $\theta$ is $K$-linear. For $\shF=\jdagY$ it is a ring
  homomorphism. Given a connection $\nabla$ on $\shF$ there is an
  obvious de Rham differential $\nabla:\Oloc^i(T,\shF)\to
  \Oloc^{i+1}(T,\shF)$ and this 
  differential is compatible with the de Rham differential on Coleman
  functions via $\theta$.
\end{proposition}
\begin{proof}
The only point which is not completely clear is the compatibility of
the differentials on $\acol$ and $\aloc$. This follows, for example
on functions, from the fact that the sections $y_x$ are horizontal,
hence for an isocrystal $M$ and a section $s\in \Hom(M,\shF)$ we
have $\nabla_{\shF}(s(y_x)) = (\nabla_{\Hom(M,\shF)}(s))(y_x)$.
\end{proof}
\begin{proposition}\label{pullinj}
  For any $z\in X$ the composition $\acol(T,\shF) \xrightarrow{z^\ast}
  \acol(T_z,\shF)
  \xrightarrow{\theta} \shF(U_z)$ is injective. In particular, the
  pullback map $z^\ast: \acol(T,\shF) \to \acol(T_z,\shF)$ is injective.
\end{proposition}
\begin{proof}
Suppose
$s(y_z)=0\in \shF(U_z)$. Let
$M_s=\iny(\Ker(s),\nabla)$. Since $y_z$ is horizontal we have by
Lemma~\ref{restriction1} that $y_z\in M_s(U_z)$. Applying
Corollary~\ref{unip-res}.\ref{unr3} to the 
inclusion $M_s\subset M$ shows that $y_x\in M_s(U_x)$ for any $x\in
X$. It follows that 
$(M_s,0,y)$ is an abstract Coleman function and the inclusion
$M_s\inject M$ defines a morphism of abstract Coleman functions. On
the other hand the $0$ map provides a map from $(M_s,0,y)$ to
$(0,0,0)$ showing that $[M,s,y]=[M_s,0,y]=0$.
\end{proof}

\begin{corollary}\label{un-prin}
  The uniqueness principle holds for Coleman functions, i.e., if
  $\theta(f)$ vanishes on an open subset of $\tu(X,P)$, then $f=0$ and
  therefore $\theta(f)$ vanishes identically.
\end{corollary}
\begin{proof}
The open subset has a non-zero open intersection with at least one
$U_x$, it follows that $\theta(f)$ vanishes on $U_x$ and by the
proposition we have $f=0$.
\end{proof}
\begin{corollary}
  The kernel of $d$ on $\acol(T)$ is $K$.
\end{corollary}
\begin{proof}
If $f\in \acol(T)$ and $df=0$ then $\theta(f)$ is locally
constant. Suppose its
value on some $U_x$ is $c$. Then $f$ and $c$ coincide on $U_x$ hence
$f=c$ by Corollary~\ref{un-prin}. 
\end{proof}
\begin{theorem}\label{exact-at-one}
  Suppose $(F,\nabla_F)$ is a unipotent isocrystal on $T$. Then the
  Coleman de Rham complex of $(F,\nabla_F)$
  is exact at the one forms. In particular, the sequence $\acol(T) \to
  \Ocol^1(T) \to  \Ocol^2(T)$ is exact.
\end{theorem}
\begin{proof}
Let $[E,\omega,y]\in \Ocol^1(T,F)$ so $\omega\in
\Hom(E,F\otimes \Omega^1)$. Define the
connection $M$ as follows. As a $\jdagY$-module
$M=E\oplus F$ and
$\nabla_M(e,f)=(\nabla_E(e),\nabla_F(f)-\omega(e))$. We let $\pi_1$
and $\pi_2$ be the projections on $E$ and $F$ respectively. Note that
$\pi_1$ is horizontal. The connection $M$
is an extension of unipotent isocrystals and is therefore
unipotent. However, it may not be integrable. We consider $N
=M^{\textup{int}}$. Choose some closed point $x_0\in X$. Since
$[E,\omega,y]$ maps to
$0$ in $\acol (T,F\otimes \Omega^2)$ it follows that
$\nabla_F(\omega(y_{x_0}))=0$.
The de Rham complex of $\nabla_F$ restricted to $U_{x_0}$ is
exact. Therefore we can choose $g\in F(U_{x_0})$ such
that $\nabla_F(g)=\omega(y_{x_0})$. Then $m_{x_0}:=(y_{x_0},g)\in
M(U_{x_0})$ and 
$\nabla_M(m_{x_0})=0$. By
Corollary~\ref{edel2cor3} we have $m_{x_0}\in N(U_{x_0})$. We can
analytically continue
$m_{x_0}$ to a collection $m=\{m_x\}$. We claim that the abstract Coleman
function $[N,\pi_2,m]$ satisfies $\nabla_F
[N,\pi_2,m]=[E,\omega,y]$. Indeed, $\nabla_F [N,\pi_2,m] =
[N,\nabla_{\Hom(M,F)}(\pi_2),m)]$ and
\begin{align*}
  (\nabla_{\Hom(M,F)}(\pi_2))(e,f)&=\nabla_F(\pi_2(e,f))-\pi_2(\nabla_M(e,f))
  \\ &= \nabla_F(f)-(\nabla_F(f)-\omega(e)) = \omega(e)
\end{align*}
so $\nabla_{\Hom(M,F)}(\pi_2)=\omega\circ \pi_1$. It follows that the
restriction of $\pi_1$ to $N$ induces a morphism $\nabla_F(N,\pi_2,m)
\to (E,\omega,y)$ and our claim follows and with it the theorem.
\end{proof}
\begin{definition}
  An abstract Coleman function $(M,s,y)$ is called minimal if the
  following two conditions are satisfied:
  \begin{enumerate}
  \item If $N\subset M$ is a subcrystal and  $y_x\in N(U_x)$ for one
    (hence all) $x\in X$, then $N=M$.
  \item There is no non-zero subcrystal of $M$ contained in $\Ker s$.
  \end{enumerate}
\end{definition}
\begin{lemma}\label{minsubquot}
  Any abstract Coleman function $(M,s,y)$ has a minimal
  subquotient.
\end{lemma}
\begin{proof}
If there is a subcrystal N such that $y_x\in N(U_x)$ then $(N,s|_N,y)$
is a subobject of $(M,s,y)$.
If $N\subset \Ker s$ then the
abstract Coleman function $(M/N,s,y \pmod{N})$ is a quotient object of
$(M,s,y)$. Repeating this we get by rank considerations to a minimal
subquotient.
\end{proof}
\begin{proposition}\label{min-unique}
  A minimal abstract Coleman function representing a given Coleman
  function is unique up to a unique isomorphism.
\end{proposition}
\begin{proof}
Suppose that $(M_i,s_i,y_i)$, $i=1,2$, represent the same Coleman
function. The Coleman function $[M_1\oplus M_2,s_1-s_2,y_1\oplus y_2]$
is therefore $0$ and it follows that for any $x\in X$, $(y_1)_x +
(y_2)_x \in \Ker (s_1-s_2)$. Let $N=\iny(\Ker(s_1-s_2),\nabla)$. Then
by Lemma~\ref{restriction1} we have for any
$x\in X$ that $(y_1)_x +(y_2)_x \in N(U_x)$. Let $\pi_i$ be the
projection from $N$ to $M_i$. We claim that both $\pi_1$ and $\pi_2$
are isomorphisms.
To show that $\pi_1$ in injective we notice that 
$\Ker \pi_1 \xrightarrow{\pi_2} M_2$ is injective and $\pi_2(\Ker
\pi_1)$ is a subcrystal of $M_2$ contained in $\Ker s_2$ hence must be
$0$ by minimality. On the other hand, $\pi_1(N)$ is a subcrystal of
$M_1$ and $(y_1)_x\in \pi_1(N)(U_x)$ for any $x\in X$ hence again by
minimality $\pi_1(N)=M_1$. The same arguments apply by symmetry to
$\pi_2$.
It is now clear that the
projections induce isomorphisms of abstract Coleman functions $\pi_i:
(N,s_i,y_1\oplus y_2) \to (M_i,s_i,y_i)$. But
$s_1=s_2$ on $N$ and therefore $(M_1,s_1,y_1)\isom (M_2,s_2,y_2)$. 
Let $\alpha_1$ and $\alpha_2$ be two maps of $(M_1,s_1,y_1)$ into
$(M_2,s_2,y_2)$. Then
$\Ker (\alpha_1-\alpha_2)$
is a subcrystal of $M_1$ containing $y_1$ and therefore must be equal to
$M_1$, hence $\alpha_1=\alpha_2$.
\end{proof}

Let $T=(X,Y,P)$ be a rigid triple. If $U\subset X$ is open, then the
triple
\begin{equation}
  \label{TU}
  T_U:=(U,Y,P)
\end{equation}
is also a rigid triple.
\begin{lemma}\label{min-res}
  Let $U$ be an open subset of $X$ and let $(M,s,y)$ be a minimal
  abstract Coleman function on $T$. Then its restriction to $T_U$ is
  also minimal.
\end{lemma}
\begin{proof}
Easy from Proposition~\ref{crist-ext}.
\end{proof}
\begin{definition}
  Let $\shF$ be a $\jdagY$-module. The association $U\mapsto
  \acol(T_U,\shF)$ defines a presheaf on the Zariski site of $X$. We
  denote this presheaf by $\Acol(T,\shF)$.
\end{definition}
\begin{proposition}
  The presheaf $\Acol(T,\shF)$ is a sheaf.
\end{proposition}
\begin{proof}
We know that for any two opens $U\subset V$ with $U$ not empty the
restriction $\acol(T_V,\shF) \to \acol(T_U,\shF)$ is injective by
Proposition~\ref{pullinj}. To prove
the sheaf property it therefore suffices to show that if $\{U_i\}$ is
an open covering of $U$ and we have $f_i\in \acol(T_{U_i},\shF)$ such
that $f_i$ and $f_j$ coincide on $T_{ij}:=T_{U_i\cap U_j}$ then there is an
$f\in \acol(T_U,\shF)$ restricting to the $f_i$. Suppose $f_i$ has a
minimal representation $(M_i,s_i,y_i)$. For each pair $i$ and $j$ the
restrictions of $(M_i,s_i,y_i)$ and $(M_j,s_j,y_j)$ to $T_{ij}$ are
minimal by Lemma~\ref{min-res} and represent the same 
function by assumption. It follows from Proposition~\ref{min-unique}
that there is a unique isomorphism $\alpha_{ij}: M_i|_{T_{ij}} \to
M_j|_{T_{ij}}$ carrying $y_i$ to $y_j$ and $s_j$ to $s_i$. By the
uniqueness these isomorphisms match well when $i$ and $j$ vary and
therefore allow us to glue the $M_i$ to an isocrystal $M$ on $T_U$.
It is
also clear than that the $s_i$ glue together to $s\in \Hom(M,\shF)$
and that the $y_i$ taken together supply a well defined horizontal
section $y_x\in M(U_x)$ for any $x\in U$. Thus we obtain our required
Coleman function.
\end{proof}
\section{Comparison with Coleman's theory}
\label{sec:comparison}

In this section we would like to show that our theory generalizes the
theory of iterated integrals due to
Coleman~\cite{Col82,Col-de88}. In Coleman's theory the functions are
built in a recursive process of integration starting from holomorphic
forms. To show that our theory gives the same result we begin by
showing that in the relevant case our theory admits a similar
recursive description. 
\begin{definition}\label{tight}
  Let $\cX\subset \cY$ be an open immersion of $\vv$-schemes such that
  $\cX$ is smooth and $\cY$ is complete. We associate to the pair
  $(\cX,\cY)$ the triple $T_{(\cX,\cY)}:=(\cX\otimes_\vv
  \kappa,\cY\otimes_\vv \kappa,\hat{\cY})$, where $\hat{\cY}$ is the
  $p$-adic completion of $\cY$. We will call such a rigid triple
  \emph{tight}. An \emph{affine} rigid triple is a tight rigid triple
  $T_{(\cX,\cY)}$ with $\cX$ affine.
\end{definition}
\begin{lemma}\label{trivializes}
  If $T$ is affine and $E$ is a unipotent isocrystal on $T$, then its
  underlying $\jdagY$-module is free.
\end{lemma}
\begin{proof}
  If $(X,Y,P)$ is affine then there is a basis of strict
  neighborhoods of $\tu(X,P)$ in $\tu(Y,P)$ which are affinoid
  (compare the proof of~\cite[Proposition 1.10]{Ber97}. The result
  therefore follows from Lemma~\ref{untriv}.
\end{proof}
We now show that in the affine case our theory admits a recursive
description.
\begin{proposition}\label{recursive}
  Suppose $T$ is affine. Let $\acoln_n(T)$ be defined recursively as
  follows: Let $\acoln_1(T)=A(T)$ and let
  $\acoln_{n+1}(T)$ be the product
  inside $\acol(T)$ of $A(T)$ with $\{f\in \acol(T),\; df\in \Omega^1(T) \cdot
  \acoln_{n-1}(T)\}$. Then $\acol(T)=\cup_n \acoln_n(T)$.
\end{proposition}
We will in fact prove a stronger result that will be needed in the
next section.
\begin{definition}
  For any locally free $\jdagY$-module $\shF$ we define the subspace
  of Coleman functions of degree at most $n$ with values in $\shF$ on $T$,
  denoted $\acoln_n(T,\shF)$ to be the subspace of all Coleman
  functions $[E,s,y]$ where $E$ has a filtration $E=F^0\supset F^1
  \supset\cdots F^n \supset F^{n+1}=0$ (i.e., of length $n+1$) by sub
  isocrystals and where
  the graded pieces are trivial connections (recall that we called a
  connection trivial if it is a direct sum of one dimensional trivial
  connections).
\end{definition}
It is clear that $\acoln_n(T,\shF)$ is an $A(T)$-submodule of
$\acol(T,\shF)$ and that $\acol(T,\shF)=\cup_n \acoln_n(T,\shF)$.
\begin{proposition}\label{recgen}
  Suppose $T$ is affine. Then $\acoln_n(T,\shF)=\acoln_n(T)\cdot
  \shF(T)$, where $\acoln_n(T)$ has been defined in
  Proposition~\ref{recursive}. In
  particular, $\acoln_n(T)$ is the same as $\acoln_n(T,\jdagY)$ as
  defined here.
\end{proposition}
\begin{proof}
Suppose $[E,s,y]$ is a Coleman function in a representation of degree
at most $n$.  Since $T$ is affine the underlying 
$\jdagY$-module of $E$ is free by Lemma~\ref{trivializes}. We may
therefore write $s=\sum g_i r_i$ with $g_i
\in \Hom(E,\jdagY)$ and $r_i \in \shF(T)$. It thus suffices to
prove the proposition for $\shF=\jdagY$, i.e., prove that
$\acoln_n(T)=\acoln_n(T,\jdagY)$. It follows from the proof of
Theorem~\ref{exact-at-one}
that a closed form in $\acoln_n(T,\Omega^1)$ has an integral in
$\acoln_{n+1}(T,\jdagY)$ and so it is clear by induction that
$\acoln_n(T)\subset \acoln_n(T,\jdagY)$. For the other direction,
suppose that $[E,s,y]$ is a Coleman function and that $E$ can be
written in a short exact sequence $0\to E_1 \to E \xrightarrow{\pi_2}
E_2\to 0$ where 
$E_1$ is trivial and $E_2$ has a
filtration as in $E$ but of length
$n$. We may further find a splitting $\pi_1: E \to E_1$, which need
not be compatible with the connection. Since $s$ can be written as
$s_1\circ \pi_1+ s_2\circ \pi_2$ with $s_i\in \Hom (E_i,\jdagY)$ it
suffices to prove the result with $s=s_i\circ \pi_i$ for $i=1,2$. We
first notice that $[E,s_2\circ \pi_2,y]=[E,s_2, \pi_2(y)]$ (since
$\pi_2$ is horizontal) so it even belongs to $\acoln_{n-1}(T)$ by the
induction hypothesis. Suppose that $s_1\in
\Hom_\nabla(E_1,\one)$. Then $\nabla^\ast (s_1\circ \pi_1)$ vanishes
on $E_1$ and therefore equals $\omega_2\circ \pi_2$ for some
$\omega_2\in \Hom(E_2,\jdagOmY^1)$. It follows that $d [E,s_1\circ \pi_1,y]=
[E,\nabla^\ast (s_1\circ \pi_1),y]=[E_2,\omega,\pi_2(y)]$. This
belongs to $\acoln_{n-1}(T,\jdagOmY^1)$ hence by the induction
hypothesis to $\Omega^1(T) \cdot \acoln_{n-1}(T)$. Note that in the
$n=0$ case we have $E_2=0$ and the same argument implies simply that
$d [E,s_1\circ \pi_1,y]=0$ hence that $[E,s_1\circ \pi_1,y]$ is a
constant. Since $E_1$ is trivial a general $s_1$ can be written as a
combination of $s_1$'s of the type that was already considered with
coefficients in $A(T)$. This completes the induction step and also
shows the case $n=0$.
\end{proof}
\begin{remark}\label{MB}
It is easily seen that a unipotent isocrystal $E$ with underlying
free $\jdagY$-module is isomorphic to an isocrystal of the form
$M_B$ defined as follows: The underlying
module is $(\jdagY)^n$. The connection depends on the $n\times n$ upper
triangular matrix $B$ with entries in $\Omega^1(T)$ and diagonal
entries $0$ and
is given by $\nabla(x_1,\ldots,x_n)=(dx_1,\ldots dx_n)+
(x_1,\ldots,x_n)\cdot B$. For such an isocrystal, having a filtration
as above of length $n+1$ means that $B$ is block upper triangular with
$n+1$ blocks (the blocks are zero matrices).
\end{remark}

We now briefly recall Coleman's theory, using the description given
in~\cite[Section~\ref{sec:integ}]{Bes98b} and some of the notation
that was established before. We first make the
general remark that Coleman allows $K$
to be any complete subfield of $\C_p$ and there is no assumption that
$K$ is discretely valued. The only essential difficulty in extending
the general theory to this case is that the Frobenius behavior of
rigid cohomology (indeed, even finite dimensionality) is not known in
this case. This is the same limitation as in Coleman's theory. In
cases where one knows this behavior we believe the other details can
be extended, though we have not checked that in detail.

We consider a pair $(\cX,\cY)$ as in Definition~\ref{tight} where
$\cY$ is assumed in addition to be a smooth 
projective and surjective scheme of relative
dimension $1$ over $\vv$. Let $T=(X,Y,P)$ be the associated rigid triple as
above. Then $Y-X=\{e_1,\ldots e_n\}$, a
finite set of points. A ``basic wide open'' $U$, in Coleman's
terminology (see~\cite[2.1]{Col-de88}), is simply a strict
neighborhood of $\tu(X,P)$ in $\tu(Y,P)$ and its 
``underlying affinoid'' is just $\tu(X,P)$. The spaces of functions
and forms which we called $A(U)$ and $\Omega^1(U)$ in~\cite{Bes98b}
are what we called here $A(T)$ and $\Omega^1(T)$.
An endomorphism $\phi:T\to T$ lifting a Frobenius endomorphisms $X\to
X$ is what we called in loc.\ cit., following Coleman, a Frobenius
endomorphism of
$U$. A theorem of Coleman~\cite[Theorem 2.2.]{Col-de88} guarantees
such a Frobenius endomorphism always exists.

For a sufficiently small $U$ we can set theoretically decompose $U$ into
a disjoint union of opens sets $U_x$ over $x\in 
Y(\kappa)$. For $x\in
X$ these sets are the usual residue classes ones defined in
\eqref{eq:disc}, which we now call residue discs, and
are each isomorphic to the open unit disc $\{|z|<1\}$, while
for $x\in \{e_1,\ldots e_n\}$ these are the intersections of $U$ with
the usual discs and are isomorphic to an open annulus $\{r<|z|<1\}$.

The differential $d:A(U_x) \to \Omega^1(U_x)$ is surjective when $U_x$
is a disc. On the other hand, when $U_x$ is an annulus there is no
integral to $dz/z$. To integrate it one needs to introduce a
logarithm, and for this one chooses a branch of the $p$-adic logarithm
and define $\alog(U_x)$ to be $A(U_x)$ if $U_x$ is a disc and 
to be the polynomial ring in the function
$\log(z)$ over $A(U_x)$ if $U_x$ is an annulus with local parameter
$z$ (the choice of the local parameter $z$ does not matter). For
either a disc or an annulus $U_x$ with a local parameter $z$ we set
$\Olog^1(U_x):=\alog(U_x)dz$. The differential
$\alog(U_x)\to \Olog^1(U_x)$ is surjective also for the
annuli, as one easily discovers by doing integration by parts of
polynomials in logs with power series coefficients.
Then one defines locally analytic
functions and one forms on $U$ by
\begin{equation*}
  \aloc(U):=\prod_x \alog(U_x),\;
  \Oloc^1(U):=\prod_x \Olog^1(U_x)\;.
\end{equation*}
There is an obvious differential $d\colon \aloc(U)\rightarrow
\Oloc^1(U)$, which is clearly surjective.

Coleman's idea is now as follows: One constructs a certain subspace $M(U)$ of
$\aloc(U)$, containing $A(U)$, which we call the space of
\emph{Coleman functions}, and
a vector space map (integration), which we denote by 
$\int$ or by $\omega\mapsto F_\omega$, from
$W(U):=M(U)\cdot \Omega^1(U)$ (product taking place inside
$\Oloc^1(U)$) to $M(U)/K\cdot 1$.
The map $\int$ is characterized by three properties:
\begin{enumerate}
\item It is a primitive for the differential in the sense that
  $dF_\omega = \omega$.
\item It is Frobenius equivariant in the sense that $\int(\phi^\ast
  \omega) = \phi^\ast\int(\omega)$.
\item If $g\in A(U)$, then $F_{dg}=g+K$.
\end{enumerate}
The construction relies on a
simple principle: If $\int$ has already been defined on some space $W$,
and $\omega\in \Oloc^1(U)$ is such that there is a polynomial $P(t)$
with $K$-coefficients such that $P(\phi^\ast)\omega=\eta\in W$,
then the conditions on the integral force the equality 
$P(\phi^\ast)F_\omega=F_\eta+\text{Const}$.
When $P$ has no roots of unity as
roots this condition fixes $F_\omega$ up to a constant.
Starting with
$W_0(U)=dA(U)$ one finds a unique way of integrating all
$\omega\in W_1(U)=\Omega^1(U)$.
One defines recursively $M_{i+1}(U):= A(U)\cdot \int(W_i(U))$ and
$W_{i+1}(U):=M_{i+1}(U)\cdot\Omega^1(U)= (\int(W_i(U)))\cdot 
\Omega^1(U)$ and checks that the principle above permits extending $\int$
uniquely to $W_{i+1}(U)$.
Finally one sets $M(U)=\bigcup_i M_i(U)$. Then clearly $W(U)=\bigcup_i W_i(U)$.
The entire theory turns out to be independent of the choice of $\phi$.
\begin{theorem}
  In the situation above there exists a ring isomorphism $\tilde{\theta}:
  \acol(T) \to M(U)$ and an isomorphism $\tilde{\theta}:\Ocol^1(T)\to
  W(U)$ compatible with the
  isomorphism on functions such that for $f\in
  \acol(T)$ or $f\in \Ocol^1(T)$ we have
  $\tilde{\theta}(f)|_{\tu(X,P)}=\theta(f)$. These isomorphisms are
  also compatible with the differential.
\end{theorem}
\begin{proof}
To define $\tilde{\theta}$ we would like to extend the definition of
the fiber functors $\omega_x$ to $x\in Y-X$. This can be done as
follows: Suppose $N\in \un(T)$. Then it is $\jdag \tilde{N}$ for
some $\tilde{N}$ defined on some strict neighborhood of $\tu(X,P)$ in
$\tu(Y,P)$. For each $x\in Y-X$ this neighborhood contains an annulus
around $x$ isomorphic to $A_r=\{r< |z| <1\}$ for some $r$. Over $A_r$
the underlying module to $\tilde{N}$ trivializes by
Lemma~\ref{untriv}. Thus it is isomorphic to some $M_B$ as in the
proof of Proposition~\ref{recursive}. Since $d:\alog(U_x) \to
\Olog^1(U_x)$ is surjective it follows easily that $\nabla_N$ has a
full set of solutions on $\tilde{N}(U_x)\otimes_{A(U_x)} \alog(U_x)$ and
the functor that sends $N$ to the set of solutions is the required
fiber functor. Therefore, given a Coleman function $f=[N,s,y]$, the
collection $y$ extends to give $y_x\in \tilde{N}(U_x)\otimes_{A(U_x)}
\alog(U_x)$ and applying the section $s$ we get an element in
$\alog(U_x)$.
Thus, the definition of $\theta$ extends immediately to provide the ring
homomorphism $\tilde{\theta}$ extending $\theta$ into $\aloc(U)$ and a
similar map on differential forms. These maps are injective because they
extends $\theta$ which is already
injective by Proposition~\ref{pullinj}. The maps $\tilde{\theta}$ on
functions and forms are clearly compatible with differentials. It
therefore suffices to prove that the 
image of $\tilde{\theta}$ is exactly $M(U)$ (a similar argument
applies to differential forms). We claim that in fact
$\tilde{\theta}(\acoln_n(T))=M_n(U)$. By the way these spaces are
constructed it is easy to see that using induction we only need to
prove the following claim: If
$\tilde{\theta}(\acoln_{n-1}(T))=M_{n-1}(U)$ and $f\in \acoln_n(T)$,
then $\tilde{\theta}(f)=\int\tilde{\theta}(df)+\textup{Const}$. This
follows by the fact that $\int$ is uniquely determined by the
condition of being $\phi$-equivariant and the fact that
$\tilde{\theta}(df) \mapsto \tilde{\theta}(f)+\textup{Const}$ is
$\phi$-equivariant since $\phi$ is an endomorphism of $T$.
\end{proof}
\section{The $p$-adic $\ddb$}
\label{sec:laplacian}

In this section we present a construction of what we call the \emph{$p$-adic
$\ddb$ operator} on Coleman forms ``of first order''. This construction
has its roots in a result of Coleman and de Shalit to be recalled
below. The justification for the title of $\ddb$ mainly comes from
the application to $p$-adic Arakelov theory~\cite{Bes00}.

Let $U$ be a basic wide open in a curve over $\C_p$. Then Coleman and de
Shalit~\cite[Lemma 2.4.4]{Col-de88} prove the following result.
\begin{proposition}
  If $\omega_i, \eta_i\in \Omega^1(U)$, $i=1,\ldots,n$, $\omega\in
  \Omega^1(U)$ and the
  $\eta_i$ are independent in $\hdr^1(U)$, then a relation
  $\sum \omega_i\cdot \int \eta_i +\omega=0$ implies that
  $\omega_i=\omega=0$ for all $i$.
\end{proposition}
We remark that the result of Coleman and de Shalit is more general
than the one we presented because it deals with integration in an
arbitrary logarithmic crystal.
\begin{corollary}
  There exists a well defined map $\ddb:W_2(U)\to \hdr^1(U)\otimes
  \Omega^1(U)$ sending $\Theta=\sum \omega_i\cdot \int \eta_i$
  to $\ddb(\Theta)=\sum[\eta_i]\otimes \omega_i$.
\end{corollary}

We call $\ddb$ the $p$-adic $\ddb$ operator. We want to
generalize the construction in terms of our new definition of Coleman
functions. As is obvious from the corollary,
our $\ddb$ should be defined on Coleman forms of degree at most
$1$. As it turns out, there is no reason to restrict to forms.

Let $\shF$ be a locally free $\jdagY$-module on $T$. A Coleman
function of degree at most $1$ on $T$ with values in
$\shF$ is given in some representation by the following data
which we encapsulate in the triple $(\short,s,y)$:
\begin{enumerate}
\item An isocrystal $E$ sitting in the short exact sequence $\short$:
  \begin{equation*}
    0\to E_1 \to E \to E_2\to 0\;,
  \end{equation*}
  such that $E_1$ and $E_2$ are trivial.
\item A homomorphism $s\in \Hom(E,\shF)$.
\item A compatible system of horizontal section $y_x\in E(U_x)$ for
  any $x\in X$.
\end{enumerate}
We perform the following construction:
The projection of the $y_x$ give a compatible system of horizontal
sections of $E_2$. Since $E_2$ is trivial
this system comes from a
global horizontal section $y_2$ of $E_2$.
The isocrystal $E$ gives an extension class $[E]\in
\Ext_\nabla^1(E_2,E_1)$. The horizontal section $y$ is an element of
$\Hom_\nabla(\jdagY,E_2)$. We can pullback the extension $[E]$ via $y_2$
to obtain $[E]\circ y_2\in
\Ext_\nabla^1(\jdagY,E_1)$. The homomorphism $s$ restricts to $s_1 \in
\Hom(E_1,\shF)$ Since $E_1$ is trivial the
natural map
\begin{equation*}
  \Hom_\nabla(E_1,\jdagY)\otimes \Hom(\jdagY,\shF)\to \Hom(E_1,\shF)
\end{equation*}
is an isomorphism. Thus we may view $s_1$ as an element of the
left hand side. There is a product
$\Hom_\nabla(E_1,\jdagY) \otimes \Ext_\nabla^1(\jdagY,E_1) \to
\Ext_\nabla^1(\jdagY,\jdagY)$. Taking this product in the first coordinate
of $s_1$ with $[E]\circ y_2$ we obtain
\begin{equation*}
  \ddb (\short,s,y) :=
  ([E]\circ y)\circ s^\prime \in \Ext_\nabla^1(\jdagY,\jdagY) \otimes
  \Hom(\jdagY,\shF)=\hr^1(X) \otimes \shF(T)\;.
\end{equation*}
A map $(\short,s,y) \to (\short\pr,s\pr,y\pr)$ is a map of abstract
Coleman function $(E,s,y)\to (E\pr,s\pr,y\pr)$ such that the map $E\to
E\pr$ sits in a map of short exact sequences $\short\to \short\pr$
\begin{lemma}\label{shortfunc}
  If there is a map $(\short,s,y) \to (\short\pr,s\pr,y\pr)$, then
  $\ddb(\short,s,y) = \ddb(\short\pr,s\pr,y\pr)$.
\end{lemma}
\begin{proof}
Formal check.
\end{proof}
\begin{proposition}
  The value of $\ddb(\short,s,y)$ depends only on the underlying
  Coleman function $[E,s,y]$ and the resulting map $\ddb:
  \acoln_1(T,\shF) \to
  \hr^1(X) \otimes \shF(T)$ is linear.
\end{proposition}
\begin{proof}
First we claim that $\ddb(\short,s,y)$ depends only on $(E,s,y)$.
Indeed, if we are given two short exact sequences for the same $E$, 
$\short$: $0\to E_1 \to E \to E_2 \to 0$ and $\short\pr$: $0\to E_1\pr
\to E \to E_2\pr \to 0$, we may consider $F$ which is the limit of the
diagram $\one \xrightarrow{(y_2,y_2\pr)} E_2\oplus E_2\pr \leftarrow
E$. There is then a map $F\to E$ which extends to a map of abstract
Coleman functions and a map of short exact sequences with both
$\short$ and $\short\pr$. Lemma~\ref{shortfunc} therefore proves the
claim. We have seen (Lemma~\ref{minsubquot}) that an abstract Coleman
function has
a minimal subquotient. Since subquotients extend to maps of short
exact sequences for an appropriate choice of filtration on the
subquotient we now see that an abstract Coleman function has the same
$\ddb$ as its minimal representative, hence $\ddb$ depends only on the
underlying Coleman function. The linearity is now straightforward.
\end{proof}
We now show that the definition of $\ddb$ above indeed coincides (up
to sign) with
the one we have given for curves using the work of Coleman and de
Shalit.
\begin{proposition}\label{affineddb}
  For an affine $T$ the operator $\ddb$ defined in the previous
  proposition sends $F\cdot f$ with $F\in \acoln_1(T)$, $dF\in
  \Omega^1(T)$ and
  $f\in \shF(T)$ to $-[dF]\otimes f \in \hr^1(X)\otimes
  \shF(T)$.
\end{proposition}
\begin{proof}
We can obtain the function $F\cdot f$ as follows: The isocrystal
is $E=M_B$ (as defined in Remark~\ref{MB}) with
$B=\left(\begin{smallmatrix} 0 & - dF\\ 0 &
  0\end{smallmatrix}\right)$, $y=(1,F)$ and $s(z,w)=w\cdot f$. We use
the obvious short exact sequence $0\to \{(0,\ast)\} \to E \to \jdagY
\to 0$. Then $y_2=1$. The pullback $[E]\circ y_2$ is obtained by finding
a preimage in $E$ to $y_2$ and applying to it the connection. Taking
the preimage $(1,0)$ we get that $E\circ y_2$ is the de Rham
cohomology class of $(0,-dF)$. The section $s_1$ is clearly given by
$\pi (f)$ with $\pi$ the isomorphism from $E_1=\{(0,\ast)\}$ to
$\jdagY$ given by the projection on the second component. The result
is now clear.
\end{proof}
\begin{proposition}\label{ddbsurjects}
  When $T$ is tight we have $\Ker \ddb = \shF(T)$. When $T$ is
  affine we have
  a short exact sequence
  \begin{equation}
    \label{eq:ddbshort}
    0\to \shF(T) \to \acoln_1(T,\shF) \xrightarrow{\ddb}
    \hr^1(X)\otimes \shF(T)\to 0
  \end{equation}
\end{proposition}
\begin{proof}
By covering a tight situation with affine ones it suffices to prove
the second statement. Suppose $T$ is affine. The surjectivity of
$\ddb$ is clear from Proposition~\ref{affineddb}. To prove the exactness
suppose 
$f\in \acoln_1^1(T,\shF)$ is in the kernel of $\ddb$. It follows from
Proposition~\ref{recgen} that $f$ can be
written as $\sum f_i \int \eta_i$ and by Proposition~\ref{affineddb} we
have then $\ddb f =
-\sum f_i \otimes [\eta_i]$. We may assume that the $f_i$
are independent over $K$. It now follows that $[\eta_i]=0$ for each
$i$, hence that $\eta_i = dg_i$ with $g_i \in A(T)$ and therefore
$\int \eta_i = g_i$ for an appropriate choice of $g_i$. Thus,
$f=\sum f_i g_i \in \shF(T)$.
\end{proof}

Consider a rigid triple $T=(X,Y,P)$. Recall that for an open $U\subset
X$ we have a rigid triple \eqref{TU} $T_U=(U,Y,P)$.
Let $\hF(T,\shF)$ be the presheaf on the Zariski site of $X$ given by
$U\mapsto \hr^1(U)\otimes \shF(T_U)$. It turns out that $\hF(T,\shF)$
is not a sheaf and that there is an interesting obstruction for
gluing. Suppose $T$ is a tight rigid triple arising from a pair
$(\cX,\cY)$. Let $\{V_i\}$ be an affine covering of $\cX$. Letting
$U_i=V_i\otimes_{\vv} \kappa$ we obtain a covering $\cU=\{U_i\}$ of
$X$. Let $\hcheck(\cU,\bullet)$ be the \ceck\ cohomology with
respect to the covering $\cU$. We shorthand $\hcheck(\cU,\shF)$ for
the \ceck\ cohomology of the presheaf $U\to \shF(T_U)$.
\begin{definition}\label{ninv}
  The map $\ninv: \hcheck^0(\cU,\hF(T,\shF))\to
  \hcheck^1(\cU,\shF)$ is the natural map derived from the exact
  sequences of \eqref{eq:ddbshort}
  \begin{equation*}
    0\to \shF(T_{U_i}) \to \acoln_1(T_{U_i},\shF) \xrightarrow{\ddb}
    \hF(T,\shF)(U_i) \to 0\;.
  \end{equation*}
\end{definition}
Explicitly, the map $\ninv$ is given as follows: For a collection
$\{\Omega_i \in \hF(T,\shF)(U_i)\}$ with $\Omega_i=\Omega_j$ on
$U_{ij}:=U_i\cap U_j$, we choose $f_i\in \acoln_1(T_{U_i},\shF)$ with
$\ddb(f_i)=\Omega_i$. The differences
$f_{ij}=f_i-f_j|_{T_{U_{ij}}}$ satisfy
$\ddb(f_{ij})=0$ hence by Proposition~\ref{ddbsurjects} $f_{ij}\in
\shF(T_{U_{ij}})$. It is immediate to check that the collection
$\{f_{ij}\}$ is a cocycle and gives a well defined cohomology class in
$\hcheck^1(\cU,\shF)$ which is $\ninv(\{\Omega_i\})$.

By \cite[1.2(ii) and Proposition 1.10]{Ber97} we have an isomorphism
between $\hr^i(X)$ and the hyper \ceck\
cohomology of $\cU$ with values in the presheaf of complexes $U\to
\Omega^\bullet(T_U)$. In particular we obtain a map $\hr^1(X) \to
\hcheck^1(\cU,A)$, where $A$ stands for the presheaf $U\to A(T_U)$.
\begin{proposition}
  We have the following commutative diagram.
  \begin{equation*}
  \xymatrix{
    {}\hr^1(X)\otimes \shF(T) \ar[r] \ar[d]
    &\hcheck^0(\cU,\hF(T,\shF)) \ar[d]^{\ninv}\\
    {}\hcheck^1(\cU,A)\otimes \shF(T) \ar[r] & \hcheck^1(\cU,\shF)
  }
  \end{equation*}
\end{proposition}
\begin{proof}
We can take an element of $\hr^1(X)$ and represent it as a
hyper-cocycle with respect to the covering $\cU$,
$\beta=(\{\eta_i\},\{g_{ij}\})$ where $\eta_i\in \Omega^1(T_{U_i})$,
$g_{ij}\in A(T_{U_{ij}})$, the $g_{ij}$ form a one-cocycle,
$d\eta_i=0$ and we have 
$dg_{ij} = \eta_i - \eta_j$ on $T_{U_{ij}}$. If $f\in \shF(T)$,
then the image of $\beta\otimes f$ in $\hcheck^0(\cU,\hF(T,\shF))$ is
given by the collection
$\Omega_i=\eta_i\otimes f|_{T_{U_i}}$. To compute the image under
$\ninv$ we, following the procedure
described after Definition~\ref{ninv}, choose lifts $f_i=(\int
\eta_i) \cdot f\in \acoln_1(T_{U_i},\shF)$ and consider the cocycle
resulting from the differences $f_i-f_j = (\int \eta_i -
\int \eta_j)\cdot f$. But $d(g_{ij} - (\int \eta_i - \int \eta_j))=0$
so each of the functions $g_{ij} - (\int \eta_i - \int \eta_j)$ is
constant on $T_{U_{ij}}$. Since the
\ceck\ cohomology of $\cU$ with constant coefficients is trivial we
can always arrange to fix the integrals $\int \eta_i$ in such a way
that $g_{ij} = \int \eta_i - \int \eta_j$ so that the application of
$\ninv$ gives $\{g_{ij}f\}$. This is clearly the same as going along the
diagram first down and then right.
\end{proof}
\begin{corollary}
  An element $\alpha\in \hcheck^0(\cU,\hF(T,\shF))$ comes from
  $\hr^1(X)\otimes \shF(T)$ if and only if $\ninv(\alpha)$ is in the image
  of $\hr^1(X)\otimes \shF(T)$ under the composed map
  $\hr^1(X)\otimes \shF(T)\to \hcheck^1(\cU,A)\otimes \shF(T) \to
  \hcheck^1(\cU,\shF)$.
\end{corollary}
\begin{proof}
The only if part is clear. For the if part we notice that we may
modify $\alpha$ by elements from $\hr^1(X)\otimes \shF(T)$ so that we
may assume $\ninv(\alpha)=0$. But then it follows from the definition
of $\ninv$ and the fact that $\acol$ is a sheaf that $\alpha$ comes
from an element $f\in \acol(T,\shF)$ via restriction to the $T_{U_i}$
and taking $\ddb$ and therefore $\alpha$ comes from $\ddb(f) \in
\hr^1(X)\otimes \shF(T)$.
\end{proof}
\newcommand{\poincare}{Poincar\'e}
\begin{remark}
  There is an interesting example to the corollary above that will be
  used in $p$-adic Arakelov theory: Suppose $T$ comes from a pair
  $(\cX,\cY)$ and $\cX=\cY$ is a relative curve and that $\shF$ is the
  sheaf $\Omega^1$. Then the composed map $\hr^1(X)\otimes \shF(T)\to
  \hcheck^1(\cU,A)\otimes \shF(T) \to \hcheck^1(\cU,\shF)$ is easily
  seen to be the same as the cup product $\hdr^1(\cX_K/K) \otimes F^1
  \hdr^1(\cX_K/K) \to \hdr^2(\cX_K/K)$. By \poincare\ duality this map
  can be non surjective only if the genus of the generic fiber is
  $0$. Thus, in positive genus the map  $\hr^1(X)\otimes \shF(T)\to
  \hcheck^0(\cU,\hF(T,\shF))$ is surjective.
\end{remark}
\section{Application: Coleman iterated integrals}
\label{sec:appl}

As an example of Coleman functions in more than one variable, we want
to discuss Coleman iterated integrals as functions of the two ends
simultaneously. The use of this for computations of syntomic
regulators of fields is sketched in~\cite{Bes-deJ98}.

Consider a tight rigid triple $T=(X,Y,P)$ of dimension $1$ and forms
$\omega_1$, $\omega_2$, $\ldots, \omega_n\in \Omega^1(T)$. Let $S$ and
$z$ be two points in $\tu(X,P)$. The iterated integrals,
\begin{equation*}
  f_k(S,z):=\int_S^z \omega_1\circ \omega_2 \circ \cdots \circ \omega_k\;,
\end{equation*}
are defined recursively, just as in the complex case, by 
\begin{equation*}
  f_1(S,z)=\int_S^z \omega_1(t),\; f_k(S,z)= \int_S^z f_{k-1}(S,t)
  \omega_k(t)\;.
\end{equation*}
The integration at each step is Coleman integration in the variable
$z$. Here we considered the point $S$ as fixed. We now want to
consider the functions $f_k(S,z)$ as a function of $S$ as well. Let
$T\times T$ be the triple $(X\times X,Y\times Y,P\times P)$.
\begin{proposition}\label{twovars}
  The functions $f_k(S,z)$ are Coleman functions in two variables on
  $T\times T$.
\end{proposition}
This is clear for $f_1$, indeed, let $F$ be a Coleman integral of
$\omega_1$. Then $f_1(S,z) =F(z)-F(S)$ which is clearly a Coleman
function.  To continue, we want to know,
at least formally, the partial derivatives of $f_k$ with respect to
$S$ and $z$. Suppose we can write $\omega_i(t)=g_i(t) dt$. Clearly,
\begin{equation*}
  \frac{\partial }{\partial z}f_{k}(S,z) = f_{k-1}(S,z)g_k(z) =
  (\int_S^z \omega_1\circ \omega_2 \circ \cdots \circ \omega_{k-1})
  g_k(z)\;.
\end{equation*}
On the other hand, we have
\[
\frac{\partial }{\partial S}f_{k}(S,z)=-f_{k-1}(S,S)\cdot g_k(S)+\int
_{z}^{S}\left(\frac{\partial }{\partial S}f_{k-1}(S,t)\right) \omega_k(t)\;.
\]
 The first term is \( 0 \). We can therefore repeat the computation expressing
the result in terms of \( f_{k-2} \) and so on. The process ends when we get
to \(\frac{\partial }{\partial S}f_{1}(S,z)=-g_1(S)\). Therefore 
\[
\frac{\partial }{\partial S}f_{k}(S,z)=-g_1(S)\int _{z}^{S}
\omega_2\circ \omega_3\circ\cdots\circ \omega_k\;.
\]
To conclude, we would expect that
\begin{equation*}
  d f_k =\Omega_k:=
  (\int_S^z \omega_1\circ \omega_2 \circ \cdots \circ \omega_{k-1})
  \cdot \omega_k(z)-
  (\int_S^z \omega_2\circ \omega_3 \circ \cdots \circ \omega_k)\cdot
  \omega_1(S)\;.
\end{equation*}
\begin{proof}[{Proof of Proposition~\ref{twovars}}]
We prove by induction on the number of differential forms the more
precise statement saying that $\int_S^z \omega_1\circ \omega_2 \circ
\cdots \circ \omega_k$ is a Coleman function of $S$ and $z$ and that
its differential is indeed $\Omega_k$ as defined above. Suppose that
we proved this for at most $k-1$ differential forms.
By the induction hypothesis it is clear that $\Omega_k$ is a Coleman
differential form and we easily compute that
\begin{align*}
  d \Omega_k&=
   d (\int_S^z \omega_1\circ \omega_2 \circ \cdots \circ \omega_{k-1})
    \wedge \omega_k(z) - 
   d (\int_S^z \omega_2\circ \omega_3 \circ \cdots \circ \omega_k)\wedge
    \omega_1(S)\\&=
    -
   (\int_S^z \omega_2\circ \omega_2 \circ \cdots \circ \omega_{k-1})
    \cdot \omega_1(S)\wedge \omega_k(z) \\&\phantom{=}\ - 
   (\int_S^z \omega_2\circ \omega_2 \circ \cdots \circ \omega_{k-1})
    \cdot \omega_k(z) \wedge \omega_1(S)=0\;
\end{align*}
It follows that $\Omega$ can be
integrated. Let $F(S,z)$ be its integral. It is immediate to see that
the restriction of $\Omega_k$ to the diagonal $S=z$ is $0$. By
functoriality $F$ is constant on the diagonal and we may therefore
assume that $F(S,S)=0$. But then for fixed $S$ this last equality
together with the fact that $d F(S,z)=(\int_S^z \omega_1\circ
\omega_2 \circ \cdots \circ \omega_{k-1}) \cdot \omega_k(z)$ implies that
$F(S,z)=f_k(S,z)$.
\end{proof}

\end{document}